\documentclass[11pt,reqno]{amsart}

\catcode`\@=11

\long\def\@savemarbox#1#2{\global\setbox#1\vtop{\hsize\marginparwidth 
  \@parboxrestore\tiny\raggedright #2}}
\newcommand\lref[1]{\ref{#1}%
\@ifundefined{r@DisplaY #1}{}{ (#1)}}

\newcommand\fakelabel[2]{\@bsphack\if@filesw {\let\thepage\relax
   \newcommand\protect{\noexpand\noexpand\noexpand}%
\xdef\@gtempa{\write\@auxout{\string
      \newlabel{#1}{{#2}{\thepage}}}}}\@gtempa
   \if@nobreak \ifvmode\nobreak\fi\fi\fi\@esphack}

\catcode`\@=12

\def\Empty{}
\newcommand\oplabel[1]{
  \def\OpArg{#1} \ifx \OpArg\Empty {} \else
        \label{#1}
  \fi}
\newtheorem{theoremSt}{Theorem}[section]

\newtheorem{exampleSt}[theoremSt]{Example}
\newtheorem{exerciseSt}[theoremSt]{Exercise}

\newcommand\MakeStEnv[1]{
  \newenvironment{#1}[1]{
  \begin{#1St} \oplabel{##1}%
  \global\def\CrntSt{\thetheoremSt}%
}{ 
  \end{#1St} }
  \newenvironment{#1+}[1]{
  \begin{#1St} \label{##1}%
  \label{DisplaY ##1}%
  \global\def\CrntSt{\thetheoremSt}%
  \def\Labl{##1}\ifx\Labl\Empty{} \else {\em (\Labl)\,}\fi%
}{ 
  \end{#1St} }
}
\MakeStEnv{theorem}
\MakeStEnv{corollary}
\MakeStEnv{proposition}
\MakeStEnv{lemma}
\MakeStEnv{definition}
\MakeStEnv{conjecture}

\long\def\State#1#2#3{
\medskip\par\noindent
{\bf #1 } {\rm(#2)}
{\it #3}
\par\medskip
}

\long\def\state#1#2{
\medskip\par\noindent
{\bf #1 } 
{\it #2}
\par\medskip
}

\newcommand{\segment}[1]{\medskip
\begin{center} \bf #1 \end{center}
\medskip}
        
\long\def\realfig#1#2#3{
\begin{figure}[htbp]
\centerline{\psfig{figure=#2}}
\caption[#1]{#3}
\oplabel{#1}
\end{figure}}

\newlength{\saveu}

\newenvironment{pf}{%
 \begin{proof}%
}{ 
 \end{proof}
}
\newenvironment{pf*}[1]{%
 \begin{proof}[#1]%
}{ 
 \end{proof}
}

\newcommand{\finishproof}[1]{ 
  \def\FPArg{#1}
  \ifx\FPArg\Empty
        \newcommand\FPArg{\CrntSt}  \fi
  \smallbreak\noindent\makebox[\textwidth]{\hfill\fbox{\FPArg}}
  \medbreak\noindent
}

\newcommand{\bfheading}[1]{\par\smallskip\noindent {\bf #1}}

\newcommand\AAA{{\mathcal A}}

\newcommand\CC{{\mathcal C}}
\newcommand\DD{{\mathcal D}}
\newcommand\EE{{\mathcal E}}
\newcommand\FF{{\mathcal F}}
\newcommand\GG{{\mathcal G}}

\newcommand\LL{{\mathcal L}}
\newcommand\MM{{\mathcal M}}
\newcommand\NN{{\mathcal N}}

\newcommand\PP{{\mathcal P}}

\newcommand\TT{{\mathcal T}}
\newcommand\UU{{\mathcal U}}

\newcommand\PMF{{\PP\kern-2pt\MM\FF}}
\newcommand\ML{{\MM\LL}}
\newcommand\PML{{\PP\kern-2pt\MM\LL}}
\newcommand\GL{{\GG\LL}}

\newcommand\half{{\textstyle{\frac12}}}

\newcommand\ep{\epsilon}

\newcommand\hhat{\widehat}
\newcommand\Proj{{\mathbf P}}

\newcommand\union{\cup}
\newcommand\intersect{\cap}
\newcommand\bbR{{\mathord{\text{I\kern-2pt R}}}}        
\newcommand\bbH{{\mathord{\text{I\kern-2pt H}}}}        

\newcommand\C{{\mathbb C}}

\newcommand\Z{{\mathbb Z}}
\newcommand\R{{\mathbb R}}

\newcommand\Hyp{{\mathbb H}}

\newcommand\PSL[1]{\text{PSL}_{#1}}

\newcommand\bigrightarrow[1]{\hbox to #1{\rightarrowfill}}
\newcommand\bigleftarrow[1]{\hbox to #1{\leftarrowfill}}

\newcommand\boundary{\partial}
\newcommand\semidir{\mathrel{\hbox{\vrule depth-.03ex height1.1ex\kern-0.15em$\times$}}}

\newcommand\til{\widetilde}

\newcommand{\diam}{\operatorname{diam}}

\numberwithin{equation}{section}

\input{psfig}

\newcommand{\T}{{\mathbf T}}

\newcommand{\pleat}{\operatorname{\mathbf{pleat}}}
\newcommand{\short}{\operatorname{\mathbf{short}}}

\newcommand{\collar}{\operatorname{\mathbf{collar}}}
\newcommand{\good}{\operatorname{\mathbf{good}}}
\newcommand{\dist}{\operatorname{dist}}

\newcommand\UML{\operatorname{\UU\MM\LL}}
\newcommand\EL{\operatorname{\EE\LL}}
\newcommand\Dnp{\DD_{\mathrm{np}}}
\newcommand\p{{\mathbf p}}
\newcommand\Teich{{\mathcal T}}

\MakeStEnv{claim}

\begin{document}

\title[Bounded geometry for Kleinian groups]{Bounded geometry for
    Kleinian groups }
\author{Yair N. Minsky}
\address{SUNY Stony Brook}
\date{\today}
\thanks{Partially supported by NSF grant DMS-9971596}

\renewcommand\marginpar[1]{} 

\begin{abstract}
We show that a Kleinian surface group, or hyperbolic 3-manifold
with a cusp-preserving homotopy-equivalence to a surface, has bounded
geometry if and only if there is an upper bound on an associated
collection of coefficients that depend only on its end
invariants. Bounded geometry is a positive lower bound on the lengths
of closed geodesics. 
When the surface is 
a once-punctured torus, the coefficients coincide with the continued
fraction coefficients associated to the ending laminations. 
\end{abstract}

\maketitle

\setcounter{tocdepth}{1}
\tableofcontents

\section{Introduction}
\label{intro}

Let $N$ be a hyperbolic 3-manifold homeomorphic to the interior of a
compact manifold. We say that $N$ has {\em bounded geometry} if,
outside the cusps of $N$, there is a positive lower bound on
injectivity radii. Equivalently, there is a lower bound on the length
of all closed geodesics in $N$.

The condition of bounded geometry is very helpful
in understanding some basic questions, such as classification by
end invariants (Thurston's Ending Lamination Conjecture), and
description of topological structure of the limit set.
For example, if $N$ is
homeomorphic to the interior of a manifold with incompressible
boundary and has no cusps, then bounded geometry 
provides an explicit
quasi-isometric model for $N$, which
is known to imply that $N$ is
uniquely determined by its end invariants, and gives a topological
model for the limit set of $\pi_1(N)$  acting on the
Riemann sphere
(see \cite{minsky:endinglam,minsky:slowmaps} and Klarreich
\cite{klarreich:thesis}; see also McMullen \cite{mcmullen:lctorus} for
progress in the unbounded geometry case). Bounded geometry also has
implications for 
the spectral theory of a hyperbolic 3-manifold and its
$L^2$-cohomology 
(Canary \cite{canary:laplacian}, Lott
\cite{lott:l2cohomology}).

The {\em end invariants} of $N$ are points in a certain parameter
space associated to each end, which describe the asymptotic geometry of
$N$. Conjecturally, $N$ is determined uniquely by these
invariants (see Thurston \cite{wpt:bull}). In this paper we address
the question of whether at least the
condition of bounded geometry can be detected from the end invariants,
and our main theorem is a first application of some tools
that were developed with the full conjecture in mind.
As a corollary we will obtain a small extension of the setting in
which the conjecture itself can be established.

\medskip

Let us restrict now to the case where $N$ is homeomorphic to
$S\times\R$, where $S$ is a surface of finite type. More specifically,
we will describe our manifolds as quotients of $\Hyp^3$ by injective
representations 
$$
\rho:\pi_1(S) \to \PSL 2(\C)
$$
with discrete images, which are also {\em type preserving}, that is
they map elements representing punctures of $S$ to parabolics. Call
such a representation a (marked) {\em Kleinian surface group}.
The theory of Ahlfors-Bers, Thurston and Bonahon attaches to $\rho$
two invariants $(\nu_+,\nu_-)$ lying in 
a combination of Teichm\"uller spaces and  lamination spaces
of $S$ and its subsurfaces.

We will associate to the pair $(\nu_+,\nu_-)$ a collection of positive
integers $\{d_Y(\nu_+,\nu_-)\}$, where $Y$  
runs over all isotopy classes of essential subsurfaces in $S$ (see
Sections \ref{complex defs} and  \ref{projection coefficients}).
These are analogues of the continued fraction coefficients, considered
in the setting where $S$ is a once-punctured torus, in \cite{minsky:torus}.
We will establish:

\state{Bounded Geometry Theorem}{
Let $\rho:\pi_1(S) \to \PSL 2(\C)$ be a Kleinian surface group with no
accidental parabolics, and end invariants $(\nu_+,\nu_-)$. Then $\rho$
has bounded geometry if and only if 
the coefficients $\{\pi_Y(\nu_+,\nu_-)\}$ are bounded above.

Moreover, for any $K>0$ there exists $\ep>0$, depending only on $K$
and the topological type of $S$, so that
$$
\sup_Y d_Y(\nu_+,\nu_-) < K \quad\implies\quad
\inf_\gamma\ell_\rho(\gamma) > \ep
$$
where the infimum is over elements of $\pi_1(S)$ that are {\em not
  externally short}, and the supremum is over essential subsurfaces
  for which $d_Y(\nu_+,\nu_-)$ is defined.
 Similarly given $\ep$ there exists $K$ for which
the implication is reversed.
}

As a corollary we obtain the following improvement of the
main theorem of \cite{minsky:endinglam}:
\State{Corollary 1}{Ending lamination theorem for bounded geometry}{
Let $\rho_1,\rho_2:\pi_1(S) \to \PSL 2(\C)$ be two  Kleinian surface
groups where $S$ is a closed surface. Suppose that
$\rho_1,\rho_2$ have the same end invariants, and that $\rho_1$ admits
a lower bound 
$$
\inf_{\gamma\in\pi_1(S)} \ell_{\rho_1}(\gamma) > 0.\eqno{(*)}
$$
Then $\rho_1 $ and $\rho_2$ are conjugate in $\PSL 2(\C)$. 
}
The main theorem  of \cite{minsky:endinglam} gives this conclusion if
{\em both} $\rho_1$ and $\rho_2$ are known to satisfy the bound $(*)$.
The Bounded Geometry theorem implies that $(*)$ is equivalent to a
condition depending only on the end invariants, which are common to
both representations. Thus the improvement closes a nagging loophole
in the earlier result. As a special case we obtain the following, 
which was not previously known: If the end invariants of $N$ are the stable and
unstable laminations of a pseudo-Anosov map $\varphi:S\to S$, then 
$N$ is the infinite cyclic cover of the hyperbolic mapping torus of $\varphi$.

The condition that $S$ be closed, rather than
finite type, is not essential, but  an extension to the case with
cusps would require a reworking of \cite{minsky:endinglam}.

\medskip

We can also improve a little on what is known about the Bers/Thurston conjecture
that every Kleinian surface group is a limit of quasifuchsian groups:

\State{Corollary 2}{Bers Density for bounded geometry}{
Any Kleinian surface group with bounded geometry and no 
parabolics is a limit of quasifuchsian groups.}

Given $\rho$ with bounded geometry,
Thurston's Double Limit Theorem (see Thurston \cite{wpt:II} and
Ohshika \cite{ohshika:ending-lams}) and a continuity property for the
lengths of measured laminations (see Brock \cite{brock:continuity,brock:boundaries})
can be used to find a sequence
$\rho_i$ of quasifuchsian representations
converging to a representation $\rho_\infty$ which has the same ending invariants
as $\rho$. Corollary 1 then implies that $\rho$ and $\rho_\infty$ are
conjugate, yielding Corollary 2. 

\medskip

Recent work of Rafi \cite{rafi:thesis} applies the results of this
paper to give a relationship between
bounded geometry for hyperbolic manifolds and Teichm\"uller geodesics.
A complete Teichm\"uller geodesic $g$ in the Teichm\"uller space
$\TT(S)$  is uniquely determined by its
endpoints $g_\pm$ in the projective measured lamination space
$\PML(S)$. 
We say that $g$ has bounded geometry if its projection to the Moduli
space of $S$ has compact closure. 

\State{Theorem}{Rafi}{
If $\rho$ is a 
Kleinian surface group and $g$ is a Teichm\"uller  geodesic such that 
$g_\pm$ and $\nu_\pm(\rho)$ are the same up to forgetting the
measures, then 
$\rho$ has bounded geometry if and only if $g$ does.
}

The ``only if'' direction was proved in \cite{minsky:endinglam}. Rafi
shows that the bounded geometry condition for a geodesic $g$ is equivalent
to a bound on $\sup_Y d_Y(g_+,g_-) $, so that this result then follows
from the Bounded Geometry Theorem.

Another consequence of Rafi's work and the work in this paper is the
following (see \cite{rafi:thesis} for details, and Farb-Mosher
\cite{farb-mosher:schottky} 
for related results). A Teichm\"uller geodesic $g$ determines a
canonical metric on $S\times\R$ as follows: For $t\in\R$ the metric on
$S\times\{t\}$ is the metric determined by the quadratic differential
of $g$ at the point $g(t)$. The distance between $S\times\{t\}$
and $S\times\{s\}$ is $|s-t|$, and the connection between slices is
given by  Teichm\"uller maps. This metric lifts to a metric on $\til
S\times\R$. 

\State{Theorem}{Rafi}{
Let $S$ be a closed surface and $g$ a Teichm\"uller geodesic in $\TT(S)$.
The lifted induced metric of $g$ on $\til S\times\R$ is
quasi-isometric to hyperbolic 3-space if and only if $g$ has bounded geometry.
}

\medskip

The Bounded Geometry Theorem also gives us a richer class of manifolds
known to have bounded geometry. Previously, aside from geometrically
finite manifolds for which bounded geometry is automatic, all examples
came from iteration of pseudo-Anosov or partial pseudo-Anosov mapping
classes acting on 
quasi-Fuchsian space or a Bers slice (see Thurston \cite{wpt:II},
McMullen \cite{mcmullen:renormbook} and Brock \cite{brock:thesis}).
The class of Teichm\"uller geodesics satisfying bounded geometry is
larger than this  -- in particular uncountable (this is a consequence 
of \cite{farb-mosher:schottky}). 
Thurston's Double Limit Theorem \cite{wpt:II} yields manifolds
whose end invariants correspond to the endpoints of these geodesics, 
and Rafi's theorem guarantees that in fact they have bounded geometry.

In spite of this one should note that bounded geometry is a rare
condition. In the
boundary of a Bers slice, for example, there is a topologically
generic (dense $G_\delta$) set of representations each of which has
arbitrarily short elements (see McMullen
\cite[Cor. 1.6]{mcmullen:cusps}, and Canary-Culler-Hersonsky-Shalen
\cite{canary-culler-hersonsky-shalen:cusps} for generalizations).
We hope that some of the techniques introduced in this paper,
when used more carefully, will yield information about general
geometrically infinite hyperbolic 3-manifolds as well.

\medskip

Some remarks on the technical conditions of the theorem: The condition of
no accidental parabolics means that all cusps of $N$ correspond to
cusps of $S$. Although accidental parabolics can be allowed, the
statement and the proof become more awkward, and we prefer to 
defer this case to a later paper.

The restriction of the second conclusion of the theorem
to curves that are not externally short rules out
those curves that are already short in the domain of discontinuity, if
any (see \S\ref{end invariants}). There are a bounded
number of such curves for any $\rho$, and we remark that at any rate they
are detected by the end invariants.

The theorem can be generalized in the standard ways to more
complicated manifolds, by considering their boundary subgroups. In the
interests of brevity we omit this discussion as well.

\newpage
\segment{Outline of the argument}

See also the research announcement \cite{minsky:rims} for a more
informal, and perhaps more readable, account of the argument. 
See Section \ref{prelims} for notation and definitions.

One direction of the theorem,
$$
\inf_\gamma\ell_\rho(\gamma) > \ep \quad\implies\quad
\sup_Y d_Y(\nu_+,\nu_-) < K 
$$
for $K$ depending only on $\ep$ and $S$ and $\gamma$ varying over the
non-parabolic elements of $\pi_1(S)$,
has already been established
in \cite{minsky:kgcc}. In fact a somewhat stronger statement is
proved, that for each individual subsurface $Y$ a lower bound on
$\ell_\rho(\boundary Y)$ implies an upper bound on $d_Y(\nu_-,\nu_+)$.

We will discuss the proof of the opposite direction,
$$
\sup_Y d_Y(\nu_+,\nu_-) < K \quad\implies\quad
\inf_\gamma\ell_\rho(\gamma) > \ep
$$
for $\ep$ depending only on $K$ and $S$, with the infimum over
$\gamma$ that are not externally short.

\medskip

For each $\gamma\in\pi_1(S)$ which is not externally short, we will
find a lower bound on $\ell_\rho(\gamma)$ by finding an upper bound on
the radius of its Margulis tube $\T_\gamma(\ep_0)$. The main idea of
the proof involves an interplay between combinatorial structure on
surfaces, namely elementary moves between
pants decompositions, and geometry of pleated surfaces, particularly
homotopies between them. The first of these is 
controlled by hyperbolicity of the complex of curves and related
results in \cite{masur-minsky:complex2,masur-minsky:complex1}, and the
second is controlled by Thurston's Uniform Injectivity Theorem
\cite{wpt:II} and its consequences.

We will construct a map $H:S\times[0,n]\to N_\rho$, whose structure is
determined by the end invariants $\nu_\pm$, and which  has the following
properties.

\begin{enumerate}
\item $H$ covers the Margulis tube $\T_\gamma(\ep_0)$ with degree 1.

\item
For each integer $i\in[0,n]$, $H_i = H|_{S\times\{i\}}$ is a pleated surface,
mapping a pants decomposition $P_i$ geodesically.

\item
$P_i$ and $P_{i+1}$ are related by an {\em elementary move}, and the
block $H|_{S\times[i,i+1]}$ has uniformly bounded tracks
$H(\{x\}\times[i,i+1])$, except  
in the special case that $\gamma$ is a component of $P_i$ or
$P_{i+1}$.

\item 
There is a uniform $M$, depending only on $\sup_Y d_Y(\nu_+,\nu_-)$,
so that $H$ is disjoint from $\T_\gamma(\ep_0)$ except on a segment
$S\times[j,j+M]$.
\end{enumerate}

Clearly, except for the special case of property (3), these properties
suffice to bound the diameter of $\T_\gamma(\ep_0)$. If the special
case occurs and $\gamma$ is a component of one of the $P_i$, 
then there is a solid torus in $S\times[0,n]$ that is mapped over
$\T_\gamma(\ep_0)$, and an additional argument must be given, using
the projection coefficient $d_\gamma(\nu_-,\nu_+)$ and a lemma on
shearing (Section \ref{shear}) to bound
the size of this solid torus. 

The initial and final pants decompositions $P_0,P_n$ are chosen using
the end invariants; for example in the quasifuchsian case they are
minimal-length pants decompositions on the convex hull boundary, and
in the degenerate case they are chosen sufficiently close to the
ending laminations to guarantee the covering property (1).
We connect $P_0$ to $P_n$ with a special sequence $\{P_i\}$ called a
{\em resolution sequence}, whose properties follow from
the work in Masur-Minsky
\cite{masur-minsky:complex2}, and
are described in Section
\ref{resolutions}. 

Property (3) is established in Section \ref{elem}, where 
we make use of the consequences of 
Thurston's Uniform Injectivity theorem 
(discussed in \S\ref{pleated surfaces}).

Property (4) is established using a quasiconvexity property of the
subset of the complex of curves spanned by curves of bounded
$\rho$-length (Section \ref{qconvex}). This together with the fact
that the resolution sequence $\{P_i\}$ follows a geodesic in the complex of
curves serves to bound the part of the sequence that can come near the
Margulis tube of $\gamma$. The dependence of the bound on 
$\sup_Y d_Y(\nu_-,\nu_+)$ is one of the properties of the resolution
sequence.

In addition to the notation and background material discussed in
Section \ref{prelims}, some well-known and/or straightforward
definitions and constructions in
hyperbolic geometry are also included in the appendix
\S\ref{appendix}.

\marginpar{
1. Examples of lams with bdd coeffs? Uncountable many?

2. Discussion of connection to other criteria (teich geods)
}

\section{Preliminaries and notation}
\label{prelims}

The following notation and constants will be used throughout the
paper. 
\begin{itemize}
\item $S$: a surface of genus $g$, with $p$ punctures, admitting a
  finite-area complete hyperbolic metric. (We may fix such a metric
  for reference, but its choice is not important).
\item $\DD(S)$: The space of discrete, faithful representations of
$\pi_1(S)$ into $\PSL 2(\C)$ that are {\em type-preserving}, meaning
that the image of any element representing a puncture is parabolic.
We call these ``Kleinian surface groups'' for short. 
\item $\Dnp(S)$: The subset of $\DD(S)$ consisting of
  representations without {\em accidental
  parabolics}, meaning that {\em only} elements representing punctures have
parabolic images. 
\item $N_\rho$: The quotient manifold $\Hyp^3/\rho(\pi_1(S))$ for a
  representation $\rho\in \DD(S)$. $N_\rho$ is homeomorphic to
  $S\times\R$ (Bonahon \cite{bonahon} and Thurston \cite{wpt:notes}),
  and comes equipped with a homotopy equivalence $S\to N_\rho$
  determined by $\rho$.

\item $l(\alpha)$: The length of a curve or arc $\alpha$. A subscript
  $l_\sigma$ usually denotes a metric, or sometimes an ambient space
  as in $l_N$.

\item $\ell(\alpha)$: The minimum of $l$ over the free homtopy class of
  a closed curve $\alpha$, or a homotopy class rel endpoints of an
  arc. Again a subscript denotes an ambient metric or space,
  and in addition $\ell_\rho$ for $\rho\in\DD(S)$ denotes
  the length of the
  shortest representative of $\alpha$ in the manifold $N_\rho$, or
  equivalently the translation length of any isometry in the conjugacy
  class determined by $\rho(\alpha)$.

\item $\ep_0$: A Margulis constant for $\Hyp^2$ and $\Hyp^3$, chosen
  as in \S\ref{margulis tubes}.

\item $\ep_1$: A number in $(0,\ep_0)$ satisfying the conditions in
  \S\ref{margulis tubes}.

\item $L_1$: A number such that every hyperbolic structure on $S$
  admits a pants decomposition of total length at most $L_1$, and
  furthermore that for any geodesic $\alpha$ in $S$ of length at least $\ep_1$ 
   such a pants decomposition exists that intersects $\alpha$. See
   \S\ref{projection coefficients}.

\item $K_0$: the bilipschitz constant in Sullivan's theorem
  (\S\ref{end invariants}).

\end{itemize}

\subsection{End invariants and convex hulls}
\label{end invariants}

For additional discussions of end invariants in the setting of this
paper see \cite{minsky:knoxville,minsky:kgcc,minsky:rims}, as well as
Ohshika \cite{ohshika:ending-lams} and of course Thurston
\cite{wpt:notes} and Bonahon \cite{bonahon}. We will recall briefly
their relevant properties in the case of no accidental parabolics,
i.e. when $\rho\in\Dnp(S)$. 

Let $\Teich(S)$ denote the Teichm\"uller space of $S$ (see e.g.
Abikoff \cite{abikoff} or Gardiner \cite{gardiner}). Let $\GL(S)$ denote
the space of geodesic laminations on $S$ and let $\ML(S)$ be the space of
measured geodesic laminations, i.e.  geodesic laminations equipped
with transverse  invariant measures of full support
(see Bonahon
\cite{bonahon:laminations} or Casson-Bleiler
\cite{casson-bleiler}). On $\GL$ we have the topology of Hausdorff
convergence of closed subsets of $S$, and on $\ML$ there is a topology
coming from  weak-* convergence of the measures on transversals. We
will  need the fact 
that if $\lambda_i\to\lambda$ in $\ML$ then the
supports of $\lambda_i$ converge in $\GL$, after restriction to a
subsequence, to a lamination containing the support of $\lambda$.

Let $\UML(S)$ denote the ``unmeasured laminations,'' 
or the quotient space of $\ML(S)$ by the
equivalence relation that forgets measures. This is a
non-Hausdorff topological space, but it has a Hausdorff subset $\EL(S)$, 
which is the image of the ``filling'' measured laminations:
those laminations   $\mu$ with the
property that $\mu$ intersects $\lambda$ nontrivially for
any $\lambda\in\ML(S)$ whose support is
not equal to the support of $\mu$. 
(See Klarreich \cite[\S 7]{klarreich:boundary} for a proof. Note
that \cite{klarreich:boundary} uses the equivalent language of
measured foliations rather than laminations).

The invariants $\nu_+$ and $\nu_-$ for $\rho\in\Dnp(S)$ lie either in
$\Teich(S)$ or in 
$\EL(S)$. Let $Q_\rho$ denote the
union of standard cusp neighborhoods for the cusps of $N_\rho$, so
that $N_\rho\setminus Q_\rho$ has two ends, which we call $e_+$ and
$e_-$. (There is an orientation convention for deciding which is
which, that will not concern us here).

Let $C(N_\rho)$ denote the convex hull of $N_\rho$. If there is a
neighborhood of the end $e_s$ (where $s$ denotes $+$ or $-$) that is disjoint
from $C(N_\rho)$ we say $e_s$ is geometrically finite, and there is a
boundary component $\boundary_sC(N_\rho)$ which is in fact the image
of a pleated surface $S\to N$, that bounds a neighborhood of $e_s$.
In the compactification $\bar N_\rho = \Hyp^3\union
\Omega_\rho/\rho(\pi_1(S))$ 
(where $\Omega_\rho$ is the domain of discontinuity in the Riemann
sphere), there is a 
component isotopic to $\boundary_sC(N_\rho)$, which inherits a
conformal structure from $\Omega_\rho$. This conformal structure, seen
as a point in $\Teich(S)$, is $\nu_s$. 

The hyperbolic structure on $\boundary_sC(N_\rho)$  yields a point 
$\nu'_s$ in $\Teich(S)$. A theorem of Sullivan
(proof in Epstein-Marden \cite{epstein-marden}) states that
$\nu'_s$ and $\nu_s$ differ by a  uniformly bilipschitz
distortion. Let $K_0$ denote this bilipschitz constant.

If $e_s$ is not geometrically finite it is geometrically infinite, and
$\nu_s$ is a lamination in $\EL(S)$, with the following properties:
\begin{enumerate}
\item There exists a sequence of simple closed curves $\alpha_i$ in $S$,
whose geodesic representatives $\alpha_i^*$ in $N_\rho$ are eventually
contained in any 
neighborhood of $e_s$ (we say ``$\alpha_i^*$ exit the end $e_s$''), and
whose lengths are bounded. 
\item For any sequence of simple closed curves $\beta_i$ whose
geodesics $\beta_i^*$ exit the 
end $e_s$, $\beta_i \to \nu_s$ in $\UML(S)$.
\end{enumerate}

\bfheading{Externally short curves.}
Call a curve $\gamma$ in $S$ {\em externally short}, with respect to a
representation $\rho$, if it is either parabolic or if
at least one of $\nu_\pm$ is a hyperbolic structure with respect to
which $\gamma$ has 
length less than $\ep_0$.

\subsection{Complexes of arcs and curves:}
\label{complex defs}

Let $Z$ be a compact surface, possibly with boundary.
If $Z$ is not an annulus, 
define  $\AAA_0(Z)$ to be 
the set of  essential homotopy classes of 
simple closed curves or properly embedded arcs in $Z$. Here
``homotopy class'' means free homotopy for closed curves, and
homotopy rel $\boundary Z$ for arcs.
``Essential'' means the homotopy class does not contain a constant
map or a map into the boundary. 

If $Z$ is an annulus, we make the same
definition except that homotopy for arcs is rel endpoints.

If $Z$ is a non-annular surface with punctures as well as boundaries,
we make a similar definition, in which arcs are allowed to terminate
in punctures.

We can extend $\AAA_0$ to a simplicial complex
$\AAA(Z)$ by letting a $k$-simplex be any $(k+1)$-tuple
$[v_0,\ldots,v_k]$ with $v_i\in\AAA_0(Z)$ distinct and having
pairwise disjoint representatives. 

Let $\AAA_i(Z)$ denote the $i$-skeleton of $\AAA(Z)$, and let
$\CC(Z)$ denote the subcomplex spanned by vertices corresponding to
simple closed curves. This is the ``complex of curves of $Z$'',
originally introduced by Harvey \cite{harvey:boundary} 
(see Harer \cite{harer:stability, harer:cohomdim} and
Ivanov \cite{ivanov:complexes1,ivanov:complexes3,ivanov:complexes2}
for subsequent developments). 

If we put a path metric on $\AAA(Z)$ making every simplex regular
Euclidean of sidelength 1, then it is clearly quasi-isometric to its
$1$-skeleton. It is also quasi-isometric to $\CC(Z)$ except in a few
simple cases when $\CC(Z)$ has no edges. When $Z$ has no boundaries or
punctures, of course $\AAA(Z)=\CC(Z)$.
Note that if $Z$ is a torus with one
hole or a sphere with four holes then this definition does not agree
with the one we used in previous papers
(e.g. \cite{minsky:torus,minsky:kgcc}) -- 
in particular $\CC(Z)$ is 0-dimensional,
whereas it was 1-dimensional before.
However the $0$-skeletons are the same in both
definitions, and the  distance function restricted from $\AAA(Z)$
agrees up to bilipschitz distortion with the distance function of the
earlier version. 


\medskip

Note that
$\AAA_0(S)$ can identified with a subset of the geodesic lamination
space $\GL(S)$. 
Let $Y\subset S$ be  a proper essential closed subsurface (all boundary curves
are homotopically nontrivial, and $Y$ is not deformable into a
cusp). We have a ``projection map''  
$$
\pi_Y : \GL(S) \to \AAA(\hhat Y)\union\{\emptyset\}
$$
defined as follows: there is a unique
cover of $S$ corresponding to the inclusion
$\pi_1(Y)\subset\pi_1(S)$, to which we can append a boundary
using the circle at infinity
of the universal cover of $S$ to yield a surface $\hhat Y$ homeomorphic
to $Y$ (take the quotient of the compactified hyperbolic plane minus 
the limit set of $\pi_1(Y)$).
Any lamination $\lambda\in\GL(S)$
lifts to this cover. If this lift has 
leaves that are either non-peripheral closed curves or essential
arcs that terminate in the boundary
components or the cusps of $\hhat Y$, these components
determine a simplex of $\AAA(\hhat Y)$
and we can take $\pi_Y(\lambda)$ to be its barycenter.
If there are no such components then $\pi_Y(\lambda)=\emptyset$.
Note that $\pi_Y(\lambda) \ne \emptyset$ whenever $\lambda$ contains a
leaf that is either a closed non-peripheral curve in $Y$, 
or intersects $\boundary Y$ essentially.

A version of this projection also appears in Ivanov
\cite{ivanov:subgroupsbook,ivanov:rank}. 

If $\beta,\gamma\in\GL(S)$  have non-empty projections $\pi_Y$,
we denote their ``$Y$-distance'' by:

$$
d_Y(\beta,\gamma) \equiv d_{\AAA(\hhat Y)}(\pi_Y(\beta),\pi_Y(\gamma)).
$$
Note that $\AAA(\hhat Y)$ can be identified naturally with $\AAA(Y)$,
except when $Y$ is an annulus, in which case the pointwise
correspondence of the boundaries matters. 
In the annulus case $d_Y$ measures relative twisting of arcs
determined rel endpoints, and in all other cases we ignore twisting on
the boundary of $\hhat Y$. If $\alpha$ is the core curve of an annulus
$Y$ we will also write 
$$
d_\alpha = d_Y.
$$
Note that, if $Y$ is a three-holed sphere (pair of pants), $\AAA(Y)$
is a finite 
complex with diameter 1, and there is not much information to be had
from $\pi_Y$. We will usually exclude three-holed spheres when
considering the projection $\pi_Y$.

\medskip

We make a final observation that one can bound
$d_{\AAA(S)}(\beta,\gamma)$, as well as $d_Y(\beta,\gamma)$ when
defined, 
in terms of the  number  of intersections of $\beta$ and $\gamma$ (although
there is no bound in the opposite direction). See e.g. Hempel
\cite{hempel:complex}. 

\medskip

\bfheading{Elementary moves on pants decompositions.}
An elementary move on a maximal curve system $P$ is a replacement of a
component $\alpha$ of $P$ by $\alpha'$, disjoint from the rest of $P$,
so that $\alpha$ and $\alpha'$ 
are in one of the two configurations shown in Figure 
\ref{elementary move fig}.

\realfig{elementary move fig}{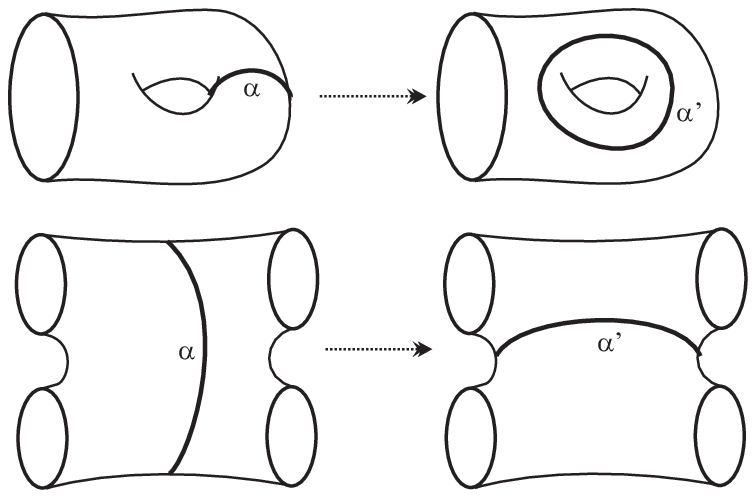}{The two types of elementary moves.}

We indicate this by  $P\to P'$ where
$P'=P\setminus\{\alpha\}\union\{\alpha'\}$ is the new curve system.
Note that there are infinitely many choices for $\alpha'$, naturally
indexed by $\Z$.

In Section \ref{elem} we will show how to relate these moves to
controlled homotopies between pleated surfaces in $N_\rho$. In Section
\ref{resolutions} we will describe the combinatorial aspects of
connecting any two pants decompositions by an efficient sequence of
elementary moves.


%

\subsection{Unmeasured laminations and the boundary of $\CC(S)$}
In addition to being
the space of ending laminations for $\rho\in\Dnp(S)$,
$\EL(S)$ also has an interpretation as the Gromov boundary
$\boundary\CC(S)$ of the
$\delta$-hyperbolic space $\CC(S)$
(see e.g.
\cite{cannon:negative,gromov:hypgroups,ghys-harpe,short:notes,bowditch:hyperbolicity} 
for material on $\delta$-hyperbolicity).

In the following theorem note that
$\CC_0(S)$ can be considered as a subset of $\UML(S)$.

\begin{theorem}{EL is boundary}
{\rm (Klarreich \cite{klarreich:boundary})}
There is a homeomorphism 
$$k:\boundary\CC(S) \to \EL(S),$$
which is natural in the sense that a sequence
$\{\beta_i\in\CC_0(S)\}$ converges to $\beta\in\boundary\CC(S)$ if and
only if it converges to $k(\beta)$ in $\UML(S)$.
\end{theorem}

\subsection{Pleated surfaces and uniform injectivity}
\label{pleated surfaces}

A {\em pleated surface} is a map $f:S\to N$ together with a hyperbolic
 metric on $S$, written $\sigma_f$ and called the {\em induced
 metric}, and a $\sigma_f$-geodesic lamination $\lambda$ on $S$,
 so that the following holds: 
 $f$ is length-preserving on paths, maps leaves of $\lambda$ to geodesics, and
 is totally geodesic on the complement  of $\lambda$. Pleated surfaces were
 introduced by Thurston \cite{wpt:notes}. See Canary-Epstein-Green
 \cite{ceg} for more details. 

The set of all pleated surfaces (in fact all maps $S\to N_\rho$)
admits a standard  equivalence relation,
in which $f\sim f\circ h$ if  $h$ is a
homeomorphism of $S$ isotopic to the identity. Let us refer to this
as equivalence {\em up to domain isotopy}.

If $P$ is a curve or arc system, i.e. a simplex in $\AAA(S)$, let
$$\pleat_\rho(P)$$
denote the set of pleated surfaces $f:S\to
N_\rho$, in the homotopy class determined by $\rho$, which map
representatives of  $P$ to geodesics. Thurston observed that such maps
always exist provided $P$ has no closed component which is parabolic
in $\rho$. Let $\pleat_\rho$ denote the set of {\em all} pleated
surfaces in the homotopy class of $\rho$.

Note, if $P$ contains an arc terminating in punctures, the
corresponding leaf in the lamination will be infinite and properly
embedded in $S$, its ends exiting the cusps.

In particular, if $P$ is a maximal curve system, or ``pants
decomposition'', $\pleat_\rho(P)$ consists of finitely many
equivalence classes, all
constructed as follows: Extend $P$ to a triangulation of $S$ with
one vertex  on each component of $P$ (and a vertex in each puncture,
if any) and ``spin'' this triangulation around $P$,
arriving at a lamination $\lambda$ whose closed leaves are $P$ and
whose other leaves spiral onto $P$, as in Figure \ref{spun
lamination}, or go out the cusps.

\realfig{spun lamination}{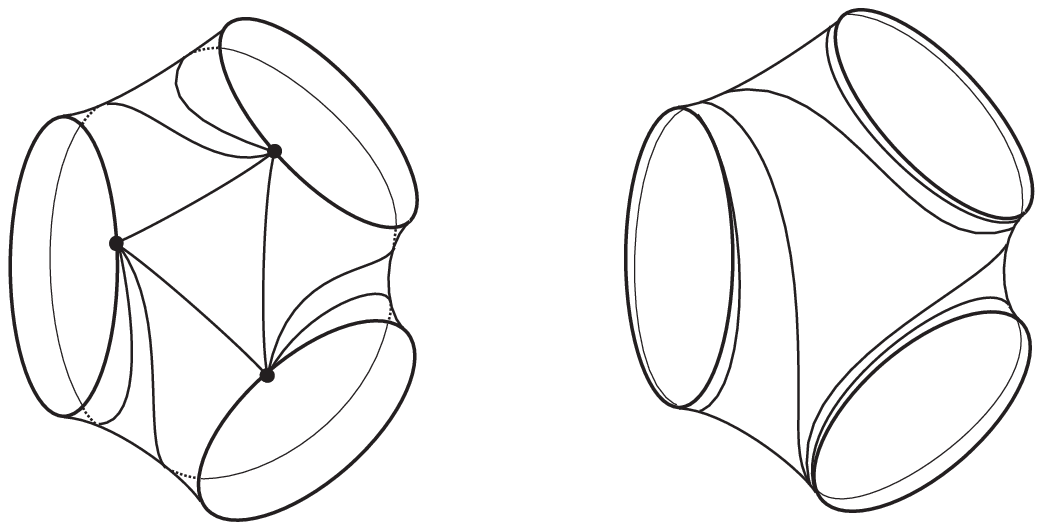}{The lamination obtained by spinning a
triangulation around a curve system. The picture shows one pair of
pants in a decomposition.}

\medskip

\segment{Uniform Injectivity}

Thurston's Uniform Injectivity theorem for pleated surfaces
\cite{wpt:I} has two corollaries that we will use here. For further
discussion and proofs see \cite{minsky:3d}, \cite{minsky:kgcc}
and also Brock \cite{brock:continuity}.

\bfheading{Bridge arcs.}
If $\alpha$ is a lamination in $S$, a {\em
  bridge arc} for $\alpha$ is an arc in $S$ with endpoints on
$\alpha$, which is not
deformable rel endpoints into $\alpha$. 
A {\em primitive bridge arc} is a bridge arc whose interior is disjoint
from $\alpha$.
If $\sigma$ is a hyperbolic metric on $S$ and $\tau$ is a bridge arc
for $\alpha$, let $[\tau]$ denote the homotopy class of $\tau$ with
endpoints fixed, and for a metric $\sigma$ let $\ell_\sigma([\tau])$
denote the length of the minimal representative of $[\tau]$.

For a lamination $\mu$ and two maps $g,g'\in\pleat_\rho(\mu)$, we say
that $g$ and $g'$ are {\em homotopic relative to $\mu$} if there is a
homotopy between them fixing $\mu$ pointwise. Lemma 3.3 in 
\cite{minsky:kgcc} guarantees that we can always precompose $g'$ by
a homeomorphism isotopic to the identity to obtain a map that is 
homotopic to $g$ relative to $\mu$.

Let $\Proj \Hyp^3$ denote the tangent line bundle over $\Hyp^3$. 
For $g\in\pleat_\rho$ mapping a lamination $\mu$ geodesically, and a
bridge arc $\tau$ of $\mu$, define $d_\p(g(\tau))$ as follows: Lift
$g(\tau)$ to an arc $\til g(\tau)$ in $\Hyp^3$ connecting two leaves
of the lift of 
$g(\mu)$. The endpoints of $\til g(\tau)$  and the leaves on which
they lie determine two points in $\Proj \Hyp^3$, and we let 
$d_\p(g(\tau))$ be their distance. In other words, $d_\p$ is small
when the two leaves are both close together and nearly tangent.

The following strengthening of Thurston's Uniform Injectivity theorem
is essentially Lemma 3.4 in
\cite{minsky:kgcc},
which follows from Lemma 2.3 in
\cite{minsky:3d}.

\begin{lemma+}{Short bridge arcs}
  Fix the surface $S$. Given $\delta_1>0$ there exists
  $\delta_0\in(0,\delta_1)$ such that the following holds. 

  Let $\rho\in\DD(S)$ and let $g\in\pleat_\rho$, 
  mapping a lamination $\mu$ geodesically.
  Suppose that $\tau$ a bridge arc for $\mu$ which is either
  primitive, or contained in the $\ep_1$-thick part of $\sigma_g$. Then
$$
d_\p(g(\tau)) \le \delta_0 \implies
\ell_{\sigma_g}([\tau]) \le \delta_1.
$$
Moreover if $g'$ is another map in $\pleat_\rho(\mu)$, 
chosen so it is homotopic to $g$ relative to $\mu$,
then
$$
\ell_{\sigma_{g}}([\tau]) \le \delta_0 \implies
\ell_{\sigma_{g'}}([\tau]) \le \delta_1.
$$
\end{lemma+}

The second statement follows from the first, since the short bridge arc
in $\sigma_{g}$ connects two leaves that contain two segments that
remain close to each other for roughly $|\log\delta_0|$, thus the
same is true for their images by $g$ (hence by $g'$). We get a
bound on $d_\p(g'(\tau))$ and up to revising the constants we have the
desired statement. 

\bfheading{Efficiency of pleated surfaces.}
Thurston also used Uniform Injectivity to establish
an estimate relating lengths of curves
in a pleated surface to their lengths in the 3-manifold.
In order to state it we need the {\em alternation number}
$$
a(\lambda,\gamma)
$$
where $\lambda $ is a lamination with finitely many leaves and
$\gamma$ is a simple closed curve (more generally a measured
lamination. This quantity, a sort of refined intersection number, 
is defined carefully in Thurston \cite{wpt:II} and Canary
\cite{canary:schottky}. For our purposes we need only the following
observations: a finite-leaved lamination consists of finitely many
closed leaves, and finitely many infinite leaves whose ends spiral
around the closed leaves.
If $\gamma$ crosses only infinite leaves of $\lambda$, then
$a(\lambda,\gamma)$ is bounded by the number of intersection points
with $\lambda$.

The following statement is a slight generalization of the theorem
proved in \cite[Thm 3.3]{wpt:II}. Thurston sketches the argument for
this generalization, and it also follows from a relative version of
the theorem proved in \cite{minsky:kgcc}.

\begin{theorem+}{Efficiency of pleated surfaces}
Given $S$ and any $\ep>0$, there is a constant $C>0$ for which the
following holds. 

Let $\rho\in\DD(S)$ and suppose $g\in\pleat_\rho$  maps geodesically
a maximal finite-leaved  lamination $\lambda$
Suppose $\gamma$ is a measured geodesic lamination in $(S,\sigma_g)$ 
which does not intersect any closed curve of $\lambda$ whose length is
less than $\ep$. Then
\begin{equation}\label{Efficiency bound}
\ell_\rho(\gamma) \le \ell_{\sigma_g}(\gamma) \le \ell_\rho(\gamma)
+ Ca(\lambda,\gamma).
\end{equation}
\end{theorem+}

\subsection{Margulis tubes}
\label{margulis tubes}

We shall denote by $\T_\alpha(\ep)$ the $\ep$-Margulis tube in
$N_\rho$ for a closed curve $\alpha$ in hyperbolic manifold, that is
the region where $\alpha$ can be represented with length at most $\ep$.
We let $\ep_0$ be a Margulis constant for 2 and 3 dimensions, meaning
that $\T_\alpha(\ep_0)$ 
is always a solid torus neighborhood of a closed geodesic, or a
horoball neighborhood of a cusp, and any two such tubes are disjoint.
We also choose $\ep_0$ sufficiently small that, on a
hyperbolic surface, any simple  closed geodesic is disjoint from any
$\ep_0$-Margulis tube but its own.

If $\rho\in\DD(S)$ is fixed, we will usually use $\T_\alpha$ to denote
$\T_{\rho(\alpha)}$, where $\alpha$ is a conjugacy class in $\pi_1(S)$.

The radius of an $\ep$-Margulis tube grows as the length of the core
curve  shrinks (see
Brooks-Matelski \cite{brooks-matelski} and Meyerhoff \cite{meyerhoff:volumes}).
We shall need to use the following facts. First, if
$\ell_\rho(\alpha)<\ep<\ep_0$ then the radius of $\T_\alpha(\ep_0)$
is at least $\half\log(\ep_0/\ep)-c$ for a universal $c$.
Second, the distance from $\T_\alpha(\ep/2)$ to
$\boundary\T_\alpha(\ep)$, if the former is non-empty, is uniformly
bounded away from $0$ and $\infty$. 

\bfheading{Margulis tubes in surface groups.}
Thurston observed that a constant $\ep_1$ exists, depending only on
$S$, so that for any $\rho\in\DD(S)$, if a pleated surface
$g\in\pleat_\rho$ meets $\T_\alpha(\ep_1)$ it can only do  so in its
own $\ep_0$-Margulis tube. Thus $\alpha$ (if it is a primitive
element) must be homotopic to the core of this tube in $S$, and in
particular simple. Together with Bonahon's tameness theorem
\cite{bonahon}, which implies that every point in $C(N_\rho)$ is within a
bounded distance of a pleated surface, we have that $\ep_1$ can be
chosen so that $\ell_\rho(\alpha) < \ep_1$ implies that $\alpha$ is
simple. 

In the remainder of the paper, we fix $\ep_1$ so that it has these
properties, and in addition $\ep_1 < \ep_0/K_0$, where $K_0$ is the
constant in Sullivan's theorem (see \S\ref{end invariants}) relating
$\nu_\pm$ to $\nu'_\pm$. This has the effect that a curve which is not
externally short has length at least $\ep_1$ in $\nu'_\pm$ as well as
$\nu_\pm$. 

\subsection{The projection coefficients}
\label{projection coefficients}

Let us now see how to define the coefficients
$$
d_Y(\nu_+,\nu_-)
$$
which appear in the main theorem, where $\nu_\pm$ are end
invariants for some $\rho\in\Dnp(S)$, and $Y$ is any essential
subsurface of $S$.

Using $\pi_Y$ as above, we can already define this 
whenever $\nu_\pm$ are laminations.
In the case of a geometrically finite end when $\nu_+$ or $\nu_-$ are
hyperbolic metrics, we can extend this definition as follows:

A theorem of Bers 
(see \cite{bers:degenerating,bers:inequality} and
Buser \cite{buser:surfaces})
says that a constant $L_1$ exists,
depending only on the topological type of $S$, so that 
for any hyperbolic metric on $S$ there is a maximal curve system (a
pants decomposition) with total length bounded by $L_1$. Moreover,
$L_1$ can be chosen so that, if $\alpha$ is a geodesic of length at
least $\ep_1$ (the constant defined in \S\ref{margulis tubes}), 
a pants decomposition can be chosen 
with length bounded by $L_1$, and intersecting  
$\alpha$ essentially. Fix this constant
for the remainder of the paper. 

Now we define
$$
\short(\sigma)
$$ 
to be the set of curve systems of $S$ with
total $\sigma$-length at most $L_1$.

Thus e.g. if both $\nu_+$ and $\nu_-$ are hyperbolic structures, we
may consider distances
$$
d_Y(P_+,P_-)
$$
for any $P_\pm\in\short(\nu_\pm)$ that both intersect $Y$ essentially,
and notice that the numbers obtained cannot vary by more than a
uniformly bounded constant (because two  different curves in 
$\short(\sigma)$ have a bounded intersection number, depending only on
$L_1$).
We let $d_Y(\nu_+,\nu_-)$ be, say,
the minimum over all choices. The case when one of
$\nu_\pm$ is a lamination and the other is a hyperbolic metric is handled
similarly. 

This defines $d_Y(\nu_+,\nu_-)$ for all $Y$,
with the exception of 
an annulus whose core curve has length less than $\ep_1$ in one of
$\nu_+$ or $\nu_-$, and a three-holed sphere all of whose boundary
curves have this property. The case of three-holed spheres will not
make any difference, since at any rate their curve complexes are finite.
The case of annuli will require a bit of attention at the end of the
proof of the main theorem.

\section{Quasiconvexity of the bounded curve set}
\label{qconvex}


Let $\CC(\rho,L)$ denote the subcomplex of $\CC(S)$ spanned by the
vertices with $\rho$-length at most $L$. In this section we will
show that this set has certain quasiconvexity properties. Before
stating them we will need to define a map 

$$
\Pi_\rho : \AAA(S) \to \PP(\CC(\rho,L_1))
$$
where $\PP(X)$ denotes the set of subsets of $X$.
Given $x\in \AAA(S)$ let $P_x$ be the curve/arc system associated to the
smallest simplex containing $x$.  We define
$$
\Pi_\rho(x) = \bigcup_{f\in\pleat_\rho(P_x)} \short(\sigma_f).
$$
For convenience we will often write $\Pi_\rho(P)$ where $P$
is a curve/arc system. 

In the following theorem, a set $A$ in a geodesic metric space is {\em
  $b$-quasiconvex} if every geodesic with endpoints in $A$ is
  contained in the $b$-neighborhood of $A$.

\begin{theorem+}{Quasiconvexity}
For any $\rho\in\DD(S)$ and $L\ge L_1$, 
$\CC(\rho,L)$ is $B$-quasiconvex, where $B$ depends only on $L$ and
the topology of $S$. 

Moreover, if $\beta$ is a geodesic in $\CC(S)$
with endpoints in $\CC(\rho,L)$ then 
$$d_\CC(x,\Pi_\rho(x)) \le B$$
for each $x\in\beta$.
\end{theorem+}

This will follow from the Coarse Projection lemma below, together with
the hyperbolicity of $\CC(S)$, via the argument in Lemma
\ref{Lipschitz implies Quasiconvex}.

\segment{The coarse projection property}

\begin{lemma+}{Coarse Projection}
For any $\rho\in\DD(S)$ the map 
$\Pi_\rho$ satisfies the following:
\begin{enumerate}
\item (Coarse Lipschitz) If $d_\AAA(x,y) \le 1$ then 
$$\diam_\CC (\Pi_\rho(x)\union\Pi_\rho(y)) \le b.$$
\item (Coarse idempotence)  If $x\in\CC(\rho,L_1)$ then
$$
d_\CC(x,\Pi_\rho(x))  = 0
$$
\end{enumerate}
where $b$ depends only on $S$.
\end{lemma+}
Note that here distance between sets (as in property (2)) is the
minimal distance, whereas $\diam_\CC$ controls maximal distances. 

\begin{pf}
Property (2) (Coarse idempotence) is immediate from the definition,
since for any $f\in\pleat_\rho(x)$, $x$ is realized in $\sigma_f$ with
its minimal length, and so is included in $\short(\sigma_f)$. 

Note next that there exists $B=B(L_1)$ such that,
for any hyperbolic structure $\sigma$, 
\begin{equation}\label{short sigma bound}
\diam_\CC(\short(\sigma)) < B
\end{equation}
since there is a uniform bound on the intersection number of any two
curves of length at most $L_1$ in the same metric on $S$. 

Now to prove part (1), it clearly suffices to consider $x,y\in\AAA_0(S)$.
The condition $d(x,y) = 1$
means $x$ and $y$ are disjoint, so $x\union y$ is a curve/arc system, and
since
$$\Pi_\rho(x\union y) \subset \Pi_\rho(x) \intersect \Pi_\rho(y),$$
the intersection is non-empty. Thus it will suffice to obtain a bound
of the form 
\begin{equation}
  \label{single Pi bound}
\diam_\CC(\Pi_\rho(\gamma)) \le B
\end{equation}
for any curve/arc system $\gamma$. 
To do this we must compare the short curves in any two different
surfaces pleated along $\gamma$. Let $f,g\in \pleat_\rho(\gamma)$, and
assume (see discussion in \S\ref{pleated surfaces}) that they are
homotopic relative to $\gamma$.

Suppose first that $\gamma$ meets the (non-cuspidal) $\ep_1$-thin part of
$\sigma_f$. Then its $f$-image meets the $\ep_1$-thin part of
$N_\rho$, and hence so does its $g$-image (since they agree). 
Let $\alpha$ be the
core curve of this component of the thin part. The length of $\alpha$
in $\sigma_g$ 
must also be at most $\ep_0$ (by the choice of $\ep_1$ in
\S\ref{margulis tubes}), and in particular
$$\alpha\in\short(\sigma_f)\intersect\short(\sigma_g).$$
Together with (\ref{short sigma bound}), this implies a bound on
$$\diam_\CC(\short(\sigma_f) \union \short(\sigma_g)).$$
The bound (\ref{single Pi bound}) follows.

Now suppose that $\gamma$ stays in the $\ep_1$-thick part of
$\sigma_f$, except possibly for cusps ($\gamma$ may contain an
infinite leaf terminating in a cusp).

By the second part of Lemma \lref{Short bridge arcs},
there exists an $\ep>0$ so that, if  
$\tau$ is a bridge arc for $\gamma$ in the $\ep_1$-thick part of
$\sigma_f$ and whose
$\sigma_f$-length is at most $\ep$, then
$\tau$ is homotopic rel endpoints to an arc of $\sigma_g$-length $\ep_0$.

Given this $\ep$, we may construct a homotopically non-trivial
curve $\gamma_\ep$ in the
$\ep_1$-thick part of $\sigma_f$, whose
$\sigma_f$-length is at 
most a constant $L_2$, and which is composed of at most two arcs of $\gamma$
and at most two primitive bridge arcs of length $\ep$ or less.
(The proof is a standard argument, which we sketch in 
Lemma \ref{short truncation}).

The bridge arcs can be homotoped to have $\sigma_g$ length at most
$\ep_0$, and hence $\gamma_\ep$ 
can be realized in $\sigma_g$ with length at most $L_2 +
2\ep_0$. In each surface this bounds its $\CC$-distance from the
curves of length $L_1$, and together with
(\ref{short sigma bound}) we again obtain a bound on 
$$\diam_\CC(\short(\sigma_f) \union \short(\sigma_g)),$$
and the desired bound (\ref{single Pi bound}) follows.
\end{pf}

\segment{Proof of Quasiconvexity}

The Quasiconvexity theorem follows from the Coarse Projection lemma
via the following standard argument, which has its roots in the proof
of Mostow's rigidity theorem (the difference between this and the
standard argument is the need to consider two projections, to a
geodesic and to the candidate quasiconvex set, whereas the standard
argument involves only projection to a geodesic).

\begin{lemma}{Lipschitz implies Quasiconvex}
Let $X$ be a $\delta$-hyperbolic geodesic metric space and $Y\subset
X$ a subset 
admitting a map $\Pi:X\to Y$ which is coarse-Lipschitz and
coarse-idempotent. That is, there exists $C>0$ such that 
\begin{itemize}
\item If $d(x,x') \le 1$ then $d(\Pi(x),\Pi(x'))\le C$, and
\item If $y\in Y$ then $d(y,\Pi(y))\le C.$
\end{itemize}
Then $Y$ is quasi-convex, and furthermore if $g$ is a geodesic in $X$  whose
endpoints are within distance $a$ of $Y$ then 
$$
d(x,\Pi(x)) \le b
$$
for some $b=b(a,\delta,C)$, and every $x\in g$.
\end{lemma}

\begin{proof}
The condition of $\delta$-hyperbolicity for $X$ implies that 
for any geodesic $g$ (finite or infinite) the closest-points
projection $\pi_g:X\to g$ 
is coarsely contracting in this sense: If $x\in X$ and $r=d(x,g)$ then
the ball $B_r(x)$ has $\pi_g$-image whose diameter is bounded by a
constant $b_0$ depending only on $\delta$. (This is an easy exercise
in the definitions -- see
e.g. \cite{gromov:hypgroups,ghys-harpe,short:notes,bowditch:hyperbolicity}.
Indeed this condition for all geodesics $g$
implies $\delta$-hyperbolicity \cite{masur-minsky:complex1}).

In the rest of the proof, let a ``dotted path'' be a sequence $p=\{x_i\in X\}$ with
$d(x_i,x_{i+1})\le C$, and let its ``length'' be $l(p) = \sum_i d(x_i,x_{i+1})$.
Let $r>C$. 
If $p$ is a dotted path in $X$ outside an $r$-neighborhood of $g$ then
the contraction property of the previous paragraph implies
\begin{equation}\label{p lower}
l(p)\ge (r-C)\left(\frac{1}{b_0}\diam(\pi_g(p)) - 1\right).
\end{equation}

Now suppose that $g$ has endpoints within $a$ of $Y$ and assume also that
$r>a$. Let $J$ be a segment of $g$ whose endpoints $\xi,\eta$ are within $r$ of
$x,y\in Y$, respectively, but whose interior is outside the
$r$-neighborhood of $Y$.  On the concatenation of geodesics
from $x$ to $\xi$, across $J$ to $\eta$ and back to $y$, select a sequence of
points spaced at most 1 apart and apply $\Pi$, obtaining (via the
coarse-Lipschitz property)
a dotted path $p$ in $Y$ satisfying
$$
l(p) \le C(2r+l(J)+1)
$$
by the coarse-Lipschitz property. 
Since $x$ and $y$ are distance $C$ from the respective  endpoints
$\Pi(x)$ and $\Pi(y)$ of $p$ (by coarse 
idempotence) and $r>C$, the contraction property of $\pi_J$ implies
that $\pi_J(\Pi(x))$ and $\pi_J(\Pi(y))$ are within $b_0$ 
of $\pi_J(x)=\xi$ and $\pi_J(y)=\eta$ respectively.
It follows that $\diam(\pi_J(p)) \ge l(J) -2b_0$.
Thus, together with (\ref{p lower}) we obtain
$$
(r-C)\left(\frac{1}{b_0}(l(J)-2b_0) - 1\right) \le
l(p) \le
C(2r+l(J)+1), 
$$
and hence
$$
\left(\frac{r-C}{b_0}-C\right)l(J) \le 3(r-C) + C(2r+1).
$$
Now if we choose $r$ so that $\frac{r-C}{b_0}-C\ge 1$,
we obtain an upper
bound
$$
l(J) \le  3(r-C)+C(2r+1).
$$
This bounds by  $b_1=r+\half(3(r-C)+C(2r+1))$ the maximum distance from
a point in $J$ to $Y$. Since this applies to every excursion of $g$
from the $r$-neighborhood of $Y$, we conclude that 
$Y$ is $b_1$-quasi-convex.

Now let $x\in g$ be any point. We have the bound $d(x,Y) \le b_1$.
Let $y\in Y$ be a nearest point to $x$. We have $d(y,\Pi(y)) \le C$ by
coarse idempotence. Now applying coarse Lipschitz to the path from $y$
to $x$, whose length is at most $b_1$, we find that
$d(\Pi(x),\Pi(y))$ is at most $C(b_1+1)$. Finally by the triangle
inequality we obtain a bound on  $d(x,\Pi(x))$.
\end{proof}

To apply this lemma to our setting, we recall first that
in \cite{masur-minsky:complex1} we proved that $\CC(S)$, and hence
$\AAA(S)$, is
$\delta$-hyperbolic. Our map $\Pi_\rho$ has images that are subsets of
$\CC(S)$ rather than single points, but this can easily be remedied by
choosing any method at all to select a single point from each set
$\Pi_\rho(x)$. Lemma \ref{Coarse Projection} implies that the
resulting map has the properties required in Lemma \ref{Lipschitz
  implies Quasiconvex}.
\qed

\section{Elementary moves on pleated surfaces}
\label{elem}

In this section we will show how to realize an elementary move between
two pants decompositions $P_0$ and $P_1$ of $S$ as a controlled 
homotopy between
pleated surfaces in $\pleat_\rho(P_0)$ and $\pleat_\rho(P_1)$. 
Lemma \lref{Homotopy bound} will show that two surfaces 
that are ``good'' with respect to a single pants decomposition admit a
controlled homotopy. 
Lemma \lref{Halfway surfaces} shows that a
pleated surface exists which is ``good'' for both $P_0$ and $P_1$
simultaneously. Thus we can concatenate a controlled
homotopy from a surface in $\pleat_\rho(P_0)$ to the halfway surface,
with one from the halfway surface to a surface in  
 $\pleat_\rho(P_1)$.

We begin with some definitions.  Let 
$\collar(\gamma,\sigma)$ denote the standard collar for $\gamma$ in
the surface $S$ with metric $\sigma$, as defined in 
Section \ref{appendix}. Similarly define $\collar(P,\sigma)$ for a
curve system $P$.

\bfheading{Good maps:}
If $P$ is a curve system in $S$,
we let 
$$\good_\rho(P,C)$$
denote the set of pleated maps $g:S\to N_\rho$
in the homotopy class determined by $\rho$, such that 
$$
 \ell_{\sigma_g}(\gamma) \le \ell_\rho(\gamma) + C
$$
for all components $\gamma$ of $P$.
Note that
$$
\pleat_\rho(P) \subset \good_\rho(P,C)
$$
for any $C\ge 0$.

\medskip

\bfheading{Good homotopies:}
Let $f,g\in\pleat_\rho$ and let $P$ be a  curve system. We say
that $f$ and $g$ admit a {\em $K$-good homotopy} with respect to $P$ if 
there exists a homotopy $H:S\times[0,1]\to N_\rho$ such that 
the following holds:
\begin{enumerate}
\item $H_0\sim f$ and  $H_1\sim g$ up to domain isotopy
(\S\ref{pleated surfaces}). 
\item Denoting by $\sigma_i$ the induced metric by $H_i$ for $i=0,1$,
$$\collar(P,\sigma_0) = \collar(P,\sigma_1).$$ 
We henceforth omit the metric when referring to these collars.
\item The metrics $\sigma_0$ and
  $\sigma_1$ are locally $K$-bilipschitz outside $\collar(P)$.
\item Let $P_0$ denote the subset of $P$ consisting of curves $\gamma$
  with $\ell_\rho(\gamma)<\ep_0$.  
The tracks $H(p\times[0,1])$ are bounded in length by $K$
  when $p\notin \collar(P_0)$.
\item For each $\alpha\in P_0$, the image
  $H(\collar(\alpha)\times[0,1])$ is contained in a $K$-neighborhood
  of the Margulis tube $\T_\alpha(\ep_0)$.
\end{enumerate}

\segment{The homotopy bound lemma}

\begin{lemma+}{Homotopy bound}
  Given $C$ there exists $K$ so that for any $\rho\in \Dnp(S)$ and
  maximal curve system $P$, if 
$$f,g \in \good_\rho(P,C)$$
 then $f $ and $g$ admit a $K$-good homotopy with respect to $P$.
\end{lemma+}
\marginpar{Extend to $\DD(S)$? Not hard}

\begin{pf}
Let us first give the proof in the case that $S$ is a closed surface. At the
end we will remark on the changes necessary to allow cusps.

Since the $\sigma_f$ and $\sigma_g$ lengths of the components of $P$
differ by at most an additive constant $C$, 
Lemma \ref{pants bilipschitz} (applied to each component of
$S\setminus P$) gives us a
homeomorphism $\varphi:S\to S$ isotopic to the identity, which takes
$\collar(P,\sigma_f) $ to $\collar(P,\sigma_g)$ and is locally
$K$-bilipschitz in its complement, with $K$ depending only on
$C$. Moreover arclengths on $\boundary\collar(P,\sigma_f)$ are {\em additively}
distorted in a bounded way: if $\alpha\subset
\boundary\collar(P,\sigma_f)$ is any arc then 
$|l_{\sigma_f}(\alpha) - l_{\sigma_g}(\varphi(\alpha))| \le K$.
After replacing $g$ with $g\circ h$, we may assume
$\collar(P,\sigma_f)=\collar(P,\sigma_g)$ (and henceforth denote it
just $\collar(P)$), and that $\sigma_f$ and $\sigma_g$ are locally
$K$-bilipschitz off $\collar(P)$, and have bounded additive length
distortion on $\boundary\collar(P)$.

Now let $H:S\times[0,1] \to N_\rho$ be the homotopy between $f$ and $g$
whose tracks $H|_{\{x\}\times[0,1]}$ are geodesics parameterized at
constant speed. ($H$ exists and is unique as a consequence of negative
curvature and the fact that $\pi_1(S)$ is non-elementary,  see
e.g. \cite{minsky:3d}).  

We will bound the tracks of $H$ on succesively larger parts of the
surface.  

\medskip

Let $Y$ denote a component of 
$S\setminus P$, and let $Y_0= Y\setminus int(\collar(\boundary Y))$. 
We first remark that, in either $\sigma_f$ or $\sigma_g$, the length
of any boundary component $\gamma$ of $Y_0$ is at most $2$ more than 
its corresponding geodesic in $S$, by (\ref{collar boundary bound}),
and this in turn is bounded by the assumption that
$f,g\in\good_\rho(P,C)$. Thus we have
\begin{equation}\label{good collar bound}
l_\sigma(\gamma) \le \ell_\rho(\gamma) + C+2
\end{equation}
for $\sigma=\sigma_f$ or $\sigma_g$.


\medskip
\bfheading{Bounds on the tripods.}
Let $\Delta$ be an
essential tripod in $Y_0$ with legs of $\sigma_f$-length bounded by $\delta$,
as given by Lemma \ref{Short tripods} in the Appendix.
Let $X$ be a 
component of the preimage of $\Delta$
in the universal cover $\til S$. 
Lift $H$ to  $\til H:\til S \times[0,1] \to \Hyp^3$.

We claim there is a uniform bound on the lengths of the tracks
$\til H(\{x\}\times[0,1])$ for $x\in X$.

\medskip

The image $\til H(X\times\{0\})$ connects three lifts of boundary
curves of $Y_0$, each invariant by a primitive deck translation 
$\gamma_i$ in $\rho(\pi_1(S))$ ($i=1,2,3$).
Let $L_i$ be the axis of $\gamma_i$ 
(see Figure \ref{Three lifts}). 
If $\gamma_i$ has translation length 
less than $\ep_1/2$, define
$N_i$ to be the lift of the corresponding $\ep_0$-Margulis tube in
$N_\rho$. Otherwise define $N_i=L_i$.

Each image (under $f$ or $g$) of a boundary curve of $Y_0$ admits a homotopy to its
geodesic representative in $N$, and because of the length bound 
(\ref{good collar bound}), Lemma \ref{Curve shortening} tells us 
that this homotopy can, in uniformly bounded distance,
be made to reach either the geodesic or its $\ep_0$-Margulis tube
if it is short. Lifting to the universal cover, we conclude that the
endpoints of $\til H(X\times\{j\})$ ($j=0,1$) are within bounded
distance of the corresponding $N_i$. 

Since we have a uniform diameter bound $d_0=2K\delta$ on $\til
H(X\times\{0\})$ and $\til H(X\times\{1\})$,
we find that 
\begin{equation}\label{triple intersection}
\til H(X\times\{j\}) \subset
\NN_{d_1}(N_1)  \intersect
\NN_{d_1}(N_2) \intersect
\NN_{d_1}(N_3) 
\end{equation}
for $j=0,1$, with a uniform $d_1$. Here $\NN_r$ denotes an
$r$-neighborhood in $\Hyp^3$.

\realfig{Three lifts}{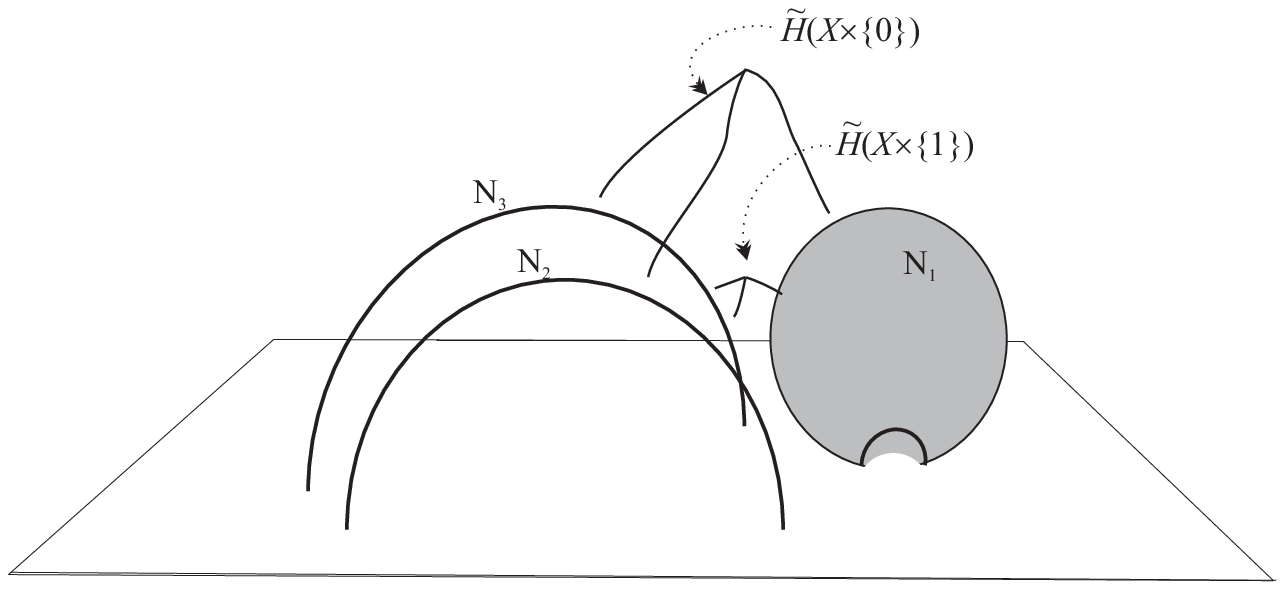}{The images of  $X\times\{0\}$
and $X\times\{1\}$  are within uniform distance of the three invariant
sets $N_1,N_2,N_3$.}

The idea now is that, by virtue of discreteness, the triple intersection 
of the  $\NN_{d_1}(N_i)$ cannot have very large diameter. Roughly, if
the $N_i$ are Margulis tubes their strict convexity implies this, and
if they are axes $L_i$ an application of Lemma \lref{Short bridge arcs} can
be used to forbid the existence of long parallel sections. Once we
prove this we will obtain a diameter bound on $\til H(X\times[0,1])$,
which will give the desired homotopy bound.


If one of the $N_i$ corresponds to a $\gamma_i$ of translation
length less than $\ep_1/2$, 
let $N_i(\ep_1)\subset N_i$ denote a lift of the corresponding
$\ep_1$-Margulis tube. 
Note that $\dist(N_i(\ep_1),\boundary N_i)$ is uniformly bounded away
from $0$ and $\infty$.
If two of the $N_i$, say $N_1$ and $N_2$, are in this case,
then they are
disjoint, and hence their convex subsets $N_i(\ep_1)$ are a definite
distance apart.
Applying Lemma \ref{convex juncture} of the Appendix to $N_1(\ep_1)$ 
and $N_2(\ep_1)$, we can deduce that
the intersection   $\NN_{d_1}(N_1)\intersect \NN_{d_1}(N_2)$ has
diameter bounded by some $d_2$.

\medskip

For the remainder of the argument, pick a pleated surface
$h\in\pleat_\rho(P)$.

Suppose that just $N_1$ corresponds to a curve of length less than $\ep_1/2$.
We claim that
$L_2$ is disjoint from $N_1(\ep_1)$. 
Let $\gamma_2^h$ denote the 
$\sigma_h$-geodesic boundary component of $Y$ corresponding to to
$\gamma_2$. If $L_2$ meets $N_1(\ep_1)$ then $\gamma_2^h$
has $h$-image that meets an
$\ep_1$-Margulis tube, and hence (by our choice of $\ep_1$)
$\gamma_2^h$ meets the $\ep_0$-collar of another boundary
component of $Y$ (in $S$ with the metric $\sigma_h$).
By our choice of $\ep_0$, a simple geodesic cannot meet an
$\ep_0$-Margulis tube in a surface unless it is the core of that tube,
so this is a contradiction.
Now $N_1(\ep_1)$ contains the convex subset
$N_1(\ep_1/2)$  which is definite distance from its boundary, 
and  hence we can apply
Lemma \ref{convex juncture} to $N_1(\ep_1/2)$ and $L_2$, deducing
a bound on
$\diam(\NN_{d_1}(N_1)\intersect \NN_{d_1}(N_2))$.

Finally suppose all three $\gamma_i$ have length at
least $\ep_1/2$. Then in fact  the
three axes $\{L_i\}$ themselves come within bounded distance $2d_3$.
If the intersection of all three $\NN_{d_3}(L_i)$ has diameter
$D$ then
$L_1,L_2$ and $L_3$ contain
segments of length at least $D-2d_3$ that remain distance $d_3$ apart.
There are therefore two a-priori constants $d_4,b>0$ so that there
exists a point $p\in\Hyp^3$
which is at most $\ep(D)=d_4e^{-bD}$ from all three $L_i$, and so that 
the tangent directions to $L_i$ at the points $x_i$ closest to $p$ are
at most $2\ep(D)$ apart in $\Proj\Hyp^3$.

Extend the tripod $\Delta$ to a tripod $\Delta'$ with endpoints in
$\boundary Y$, and let $X'$ be the component of its lift to $\til S$ that contains
$X$. The endpoints of $X'$ are mapped by $\til h$ to points $y_i\in
L_i$ ($i=1,2,3$). The arcs $[y_i,x_i]$ on $L_i$ pull back to arcs on
$\boundary Y$ which, if we append them to $\Delta'$ and perturb
slightly to an embedding in $Y$, give us a new tripod $\Delta''$ whose
endpoints (after lifting and applying $\til h$) map to $x_i$. 

Let $\delta'$ be the constant given in Lemma \ref{Short tripods}, and
let $\delta_0$ be the constant given by Lemma \ref{Short bridge arcs}
after setting $\delta_1=\delta'$. Note that each pair of legs of
$\Delta''$ is a primitive bridge arc for $\boundary Y$, whose $h$-image is
homotopic rel endpoints to an arc of length at most $2\ep(D)$.
If $D$ is sufficiently large that $2\ep(D)<\delta_0$, then Lemma
\ref{Short bridge arcs} tells us that each pair of legs of $\Delta''$
is homotopic to an arc of length at most $\delta_1$ in
$\sigma_h$. This gives us a triangle in $Y$ with the same vertices as
$\Delta''$ whose sidelengths are at most $\delta_1$. Joining its
barycenter to its vertices we obtain a new tripod $\Delta'''$ all of
whose legs are bounded by $\delta_1=\delta'$. This contradicts Lemma
\ref{Short tripods}. We conclude that the triple intersection of
(\ref{triple intersection}) has diameter at most $D$.

%

This diameter bound and (\ref{triple intersection}) now imply a
uniform bound on the track lengths of
$\til H$ restricted to $X$, and hence $H$ restricted to $\Delta$.

\medskip
\bfheading{Bounds outside the collars.}
To bound the tracks of $H$ on the rest  of $S\setminus
\collar(P)$, we first bound them on $\boundary\collar(P)$.

Let $\gamma$ be a component of $\boundary \collar(P)$. Let $Y$ be the
component of $S\setminus P$ containing $\gamma$, so that
$\gamma$ is in the boundary of $Y_0\equiv Y\setminus int(\collar(\boundary
Y))$. Lemma \ref{Short tripods} gives us a bounded essential tripod
$\Delta$ in $Y_0$ with endpoint $x\in \gamma$.
By the forgoing discussion
we already have a bound, say $t_0$, on the
track length $H(\{x\}\times[0,1])$. 

If $\ell_\rho(\gamma)\le 1$, we have the bound (\ref{good collar bound})
on the $l_{\sigma_f}(\gamma)$ and $l_{\sigma_g}(\gamma)$.
Via the triangle inequality
we immediately get a bound on $H(\{y\}\times[0,1])$ for any
$y\in\gamma$.

\medskip

Suppose that $\ell_\rho(\gamma) > 1$. Let $G = H|\gamma\times[0,1]$,
and let $\til G:\til\gamma\times[0,1] \to \Hyp^3$ be a lift to the
universal covers. Lift $x$
to $\{x_i|i\in\Z\}$ in $\til \gamma$. 
Let $\Gamma$ be the geodesic lift of $\gamma^*$
to $\Hyp^3$ which is invariant by the holonomy of $\til G$. 

Because of the bound (\ref{good collar bound}) on $l_\sigma(\gamma)$
(for $\sigma=\sigma_f$ or $\sigma_g$), we may apply 
Lemma \ref{Curve shortening} to either $\til G_0 = \til G(\cdot,0)$ or
$\til G_1= \til G(\cdot,1)$. 
For any $y\in\til\gamma$  let $\eta_i=\pi(\til G_i(y))$, where
$\pi:\Hyp^3\to \Gamma$ is the orthogonal projection. Then the lemma
gives us, first, a bound
$$d(\til G_i(y),\eta_i) \le a_1$$
and $i=0,1$. 

Let $\xi_{ij} = \pi(\til G_i(x_j))$. The main part of  Lemma \ref{Curve
shortening} tells us that
for any $y\in[x_0,x_1]$, its projection
$\eta_0$ lies a bounded distance $a_2$ from the point in $[\xi_{00},\xi_{01}]$
that is at distance $|[x_0,y]|_f$ from $\xi_{00}$. Here
 $|\cdot|_f$ and $|\cdot|_g$ denote the lifts of
$\sigma_f$ and $\sigma_g$ arclength to $\til\gamma$. Similarly we have
a bound for $\eta_1$,  $[\xi_{10},\xi_{11}]$
and $|[x_0,y]|_g$.

The length bound on $G(x\times[0,1])$ tells us that
$d(\xi_{0j},\xi_{1j})\le 2a_1+t_0$.
Since $|[x_0,y]|_f$ and $|[x_0,y]|_g$ differ by at most $K$
(via Lemma \ref{pants bilipschitz}, as discussed in the beginning of
the proof), we obtain a uniform upper bound on $d(\eta_0,\eta_1)$ and
hence on the track length $\til G(\{y\}\times[0,1])$.

This proves our uniform bound on the  track lengths for
$\boundary\collar(P)$.

\medskip

Now for any point $z$ in $S\setminus\collar(P)$ there is an arc $\beta$
connecting $z$ to 
some  point $z'$ in $\boundary\collar(P)$ whose length in
either $\sigma_f$ or $\sigma_g$ is uniformly bounded. 
The length of $H(z'\times[0,1])$ is bounded by the above, so 
we may bound $H(z\times[0,1])$ using the triangle inequality.

\medskip
\bfheading{Bounds in the collars.}
It remains to control $H$ on 
$\collar(P)\times[0,1]$. For each component $\gamma$ of $P$, 
$\collar(\gamma)\times[0,1]$ is a solid torus, and we will control $H$
by considering a meridian curve 
$$
m = \boundary (a\times[0,1])
$$
where $a$ is an arc connecting the boundary components of
$\collar(\gamma)$. 

Suppose first that $\ell_\rho(\gamma) \ge \ep_0$. Then the radius of
$\collar(\gamma)$ is bounded above by $s\equiv w(\ep_0)$ in either
$\sigma_f$ 
or $\sigma_g$. Let $a$ denote a minimal $\sigma_f$-length arc crossing
$\collar(\gamma)$, 
and let us find a bound on $L\equiv l_{\sigma_g}(a)$. (Note that up to the
usual equivalence of maps we may assume $a$ is geodesic in both
metrics). 

Let $\gamma_g$ be the $\sigma_g$-geodesic representative of $\gamma$.
We first replace $a$ by a homotopic path (with fixed endpoints)
$q_0*a'*q_1$ where $a'$ is the $\sigma_g$-orthogonal 
projection of $a$ to $\gamma_g$, and $q_i$ are orthogonal segments to
$\gamma_g$ of length at most $s$. Thus $l(a') \ge L-2s$.
Now apply Lemma \ref{Curve shortening} to deform $g(a')$, fixing
endpoints,  to a path
$p_0*a''*p_1$ in $N$ with  $p_i$ of length bounded by $a_1$,
$a''$ traveling along the geodesic representative of $\gamma$, 
$$l_N(a'') \ge c_1 l(a') - c_2
$$
where $c_1$ and $c_2$ depend on the constants in the lemma.
It follows that the geodesic representative of the path $g(a)$ in $N$
has length at least $c_3L-c_4$ for suitable constants $c_3,c_4$.

Now consider a lift $\til H(m)$ of the meridian to $\Hyp^3$. The
endpoints of $\til H(a\times\{1\})$ are at least $c_3L-c_4$ apart,
and on the other hand the other three legs of $m$
give us an upper bound of $2s + 2t_1$, where $t_1$ is the track bound we
have already obtained on $H$ outside the collars. This gives us an
upper bound on $L$, and hence on the length of $H(m)$.

A bound on the tracks of $H$ follows immediately for any point of $a$, and
since we can foliate $\collar(\gamma)$ by arcs such as $a$, we obtain
a bound in the entire collar.

\medskip

Finally suppose that $\ell_\rho(\gamma) < \ep_0$, that is $\gamma\in P_0$.
Since $\ell_{\sigma_f}(\gamma)$ and $\ell_{\sigma_g}(\gamma)$ are
bounded by $\ep_0+C$, we may foliate $\collar(\gamma)$ by closed curves with
this length bound in both metrics (again, up to precomposing the maps
by homeomorphisms homotopic to the identity, we may assume the same
curves are bounded in both metrics).
For each such curve $\beta$, the geodesic
homotopy $H$ must be contained in $\T_\gamma(\ep_0)$ for all but a
bounded portion. This is a standard area argument, for the area of the
ruled annulus $H(\beta\times[0,1])$ is bounded by the length of its
boundary, and on the other hand a long section of the annulus outside
of $\T_\gamma(\ep_0)$ would have area at least $\ep_0$ times its length.
(See \cite{wpt:notes,bonahon} for similar area arguments).

Thus  $H(\collar(\gamma)\times[0,1])$ is contained in a uniformly
bounded neighborhood of $\T_\gamma(\ep_0)$, which is what we needed to prove.

\bfheading{Surfaces with cusps.}
It remains to discuss the case when $S$ has cusps. The main difference
is that the components $Y$ of $S\setminus P$ may have cusps rather
than boundary components. All the arguments go through in the same
way, with $\collar(P)$ replaced by the union of $\collar(P)$ with the
collars of the cusps. In the complement of these collars we still
obtain a bilipschitz relation between $\sigma_g$ and $\sigma_f$.
In bounding the tracks on tripods, the
neighborhoods $N_i$ must be allowed to be Margulis tubes of parabolic
elements when the corresponding boundary component of $Y$ is a
cusp. 
\end{pf}

\bfheading{Remark:} The bilipschitz relation between $\sigma_f$ and
$\sigma_g$ generally breaks down in the collars of $P$, but a careful
consideration of the proof will show that, if $\gamma\in P$ is a
component with $\ell_\rho(\gamma)$ bounded both above and below, then
we may extend the bilipschitz relation to $\collar(\gamma)$ as
well. This is a consequence of the simultaneous bound on the  arcs $a$
crossing the collar in both metrics.

\segment{Halfway surfaces}

\begin{lemma+}{Halfway surfaces}
There exists a $C_1>0$ depending only on $S$, so that if $\rho\in\DD(S)$
and $P_0\to P_1$ is an elementary move then 
$$
\good_\rho(P_0,C_1)\intersect \good_\rho(P_1,C_1) \ne \emptyset.
$$
\end{lemma+}
A map in this intersection is called a 
{\em  halfway surface} for $P_0$ and $P_1$.

\begin{pf}
We will construct a pleated surface
$g\in \pleat_\rho(P_0\intersect P_1)$ to which we can  apply
Thurston's Efficiency of Pleated Surfaces 
(Theorem \ref{Efficiency of pleated surfaces}).

Let $\alpha_0\in P_0$ and $\alpha_1\in P_1$ be the curves exchanged by
the elementary move. 
Let $Y$ be the component of $S-(P_0\intersect P_1)$ containing
$\alpha_0$ and $\alpha_1$. Note that $Y$ is a 4-holed sphere or
1-holed torus. 

%
%
%
%

Let us describe a lamination $\lambda$ on $Y$ by first considering its
lift to a planar cover. In Figure \ref{lambda in cover}
we indicate the plane $\R^2$ with a
small disk removed around every point in the lattice $\Z^2$. 
When $Y$ is a one-holed torus it is obtained as the quotient of this
by the action of $\Z^2$, and when $Y$ is a four-holed sphere it is the
quotient by the group generated by $(2\Z)^2$ and $-I$. Normalize the
picture so that $\alpha_0$ lifts to lines parallel to the $x$-axis and
$\alpha_1$ lifts to lines parallel to the $y$-axis.

\realfig{lambda in cover}{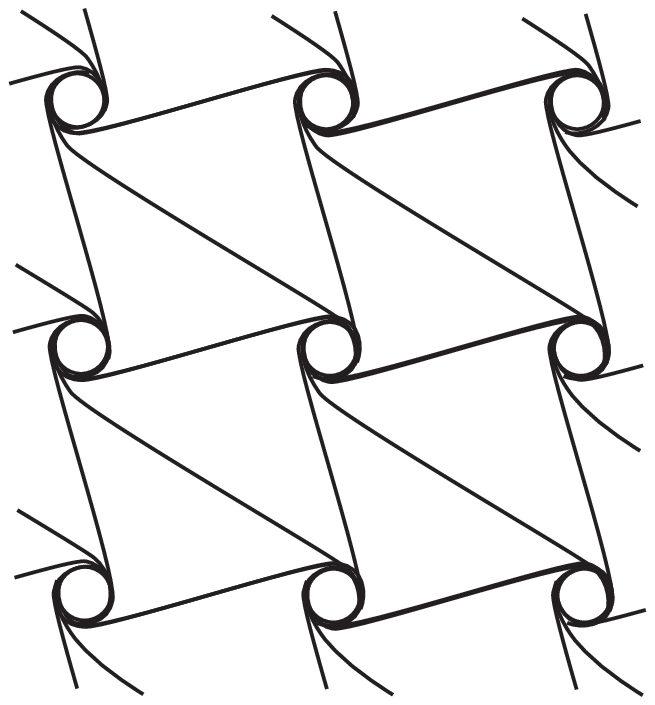}{The lamination $\til\lambda$
in the planar cover of $Y$.}

We have indicated the lamination $\til\lambda$ in this cover, which is
obtained from a standard triangulation by spinning leftward around the
boundary components. The projection of this to $Y$ is $\lambda$.
This discussion applies equally well
when $Y$ has ends that are cusps of $S$. In this case the leaves that
we drew as winding around $\boundary Y$ simply go out the
corresponding cusp.

Let $g\in \pleat_\rho(P_0\intersect P_1)$ be the pleated surface
mapping $\lambda$ geodesically.
An inspection of the diagram gives this bound on the alternation
numbers: 
$$
a(\lambda,\alpha_i) \le 4
$$
for $i=0,1$ (in fact it is 2 when $Y$ is a one-holed torus and $4$
when $Y$ is a four-holed sphere).
And furthermore that $\alpha_0$ and $\alpha_1$ cross no closed leaves
of $\lambda$. Efficiency of Pleated Surfaces then gives us the
inequality
$$
\ell_{\sigma_g}(\alpha_i) \le \ell_\rho(\alpha_i) + C 
$$
for $i=0,1$ and a uniform $C$. Thus $g\in \good_\rho(P_i,C)$ for
$i=0,1$, and the lemma is proved.
\end{pf}

%
%
%

The following lemma controls the geometry of a halfway surface. It
will be instrumental in the last part of the proof of the main
theorem.

\begin{lemma}{halfway annulus bound}
Let $P_0\to P_1$ be an elementary move exchanging $\alpha$ and 
$\alpha'$, and let
$g\in\good_\rho(P_0,C_1)\intersect\good_\rho(P_1,C_1)$.
Assume that $\alpha$ and $\alpha'$ are realized as geodesics in
$\sigma_g$, and suppose that
$$
\ell_\rho(\alpha) < \ep_1/2
$$
but 
$$
\ell_{\sigma_g}(\alpha) \ge \ep_1.
$$
Then there is an upper bound $C_2$, depending only on $\ep_1, C_1$ and
the topological type of $S$, so that 
$$
l_{\sigma_g}(a) \le C_2
$$
for each arc $a$ of $\alpha'\intersect\collar(\alpha,\sigma_g)$.
\end{lemma}

\begin{proof}

Let $a$ be an arc of 
$\alpha'\intersect\collar(\alpha,\sigma_g)$ -- there may be one or two
such arcs. Let $L=l_{\sigma_g}(a)=l_N(g(a))$. 
Deform $a$ fixing endpoints to a path $q_0*a'*q_1$
with $q_i$ orthgonal paths from $\boundary a$ to $\alpha$, and $a'$
running along $\alpha$. The lengths of $q_i$ are at most $w(\ep_0)$,
and $l(a') \le L$.
Lemma \ref{Curve shortening} allows us to deform $g(a')$, fixing
endpoints, to $p_0*a''*p_1$ with $l_N(p_i)\le a_1$ and $a''$ in the
$\ep_1/2$-Margulis tube of $\alpha$, so that
$$
l_N(a'') = n \ep_1/2  + r'
$$
where we write $l(a') = n\ell_{\sigma_g}(\alpha) + r$, $n\in\Z_{\ge
0}$, $r\in[0,\ell_{\sigma_g}(\alpha)]$, and $|r-r'|\le a_2$.
Since $\ell_{\sigma_g}(\alpha)\ge \ep_1$, we conclude that $g(a)$ can
be deformed to an arc 
$a''' = g(q_0)*p_0*a''*p_1*g(q_1)$, so that
\begin{equation}
l_N(a''') \le \half L + c 
\end{equation}
with $c$ depending on the previous constants.
This is shorter than $g(a)$ by at least $L/2 - c$, so 
we can  shorten the curve $g(\alpha')$ by at least this much,
concluding that
$$
\ell_\rho(\alpha') \le \ell_{\sigma_g}(\alpha') - \half L+c.    
$$
On the other hand $\ell_{\sigma_g}(\alpha') - C_1 \le \ell_\rho(\alpha') $
since $g\in\good_\rho(\alpha',C_1)$, so  we obtain an upper bound 
on $L$.
\end{proof}

\section{Resolution sequences} 
\label{resolutions}

In Masur-Minsky \cite{masur-minsky:complex2}, we show the existence of special
sequences of elementary moves that are controlled in terms of the
geometry of the complex of curves, and particularly the projections
$\pi_Y$. First some terminology: if $P_0\to P_1 \to \cdots \to P_n$ is
an elementary-move sequence and $\beta$ is any vertex of $\CC(S)$, we
denote
$$J_\beta = \{ j\in[0,n]: \beta\in P_j\}.$$
Note that if  $J_\beta$ is an interval $[k,l]$,
then the elementary move $P_{k-1}\to P_k$ exchanges some $\alpha$ for
$\beta$, and $P_l\to P_{l+1}$ exchanges $\beta$ for some
$\alpha'$. Both $\alpha$ and $\alpha'$ intersect $\beta$, and we call
them the {\em predecessor} and {\em successor} of $\beta$, respectively.

In the following theorem, where $\beta_0,\ldots,\beta_m$ is a sequence of
vertices in $\CC(S)$ we use the notation
$$J_{[s,t]} \equiv \bigcup_{i=s}^t J_{\beta_i}.$$

\begin{theorem+}{Controlled Resolution Sequences}
Let $P$ and $Q$ be maximal curve systems in $S$. 
There exists a geodesic in $\CC_1(S)$ with vertex sequence
$\beta_0,\ldots,\beta_m$, and an elementary move sequence
$P_0\to\ldots\to P_n$, with the following properties:
\begin{enumerate}
\item $\beta_0\in P_0 = P$ and $\beta_m \in P_n = Q$.
\item Each $P_j$ contains some $\beta_i$.
\item $J_\beta$, if nonempty, is always an interval, and if
$[s,t]\subset[0,m]$ then
$$
|J_{[s,t]}| \le K (t-s) \sup_Y d_Y(P,Q)^a,
$$
where the supremum is over only those non-annular subsurfaces $Y$
whose boundary 
curves are components of some $P_k$ with $k\in
J_{[s,t]}$. 
\item
If $\beta$ is a curve with non-empty $J_\beta$, then its predecessor
and successor curves $\alpha $ and $\alpha'$ satisfy
$$
| d_\beta(\alpha,\alpha') - d_\beta(P,Q) | \le \delta.
$$
\end{enumerate}
The constants $K,a,\delta$ depend only on the topological type of $S$.
The expression $|J|$ for an interval $J$ denotes its diameter.
\end{theorem+}
The sequence $\{P_i\}$ in this theorem is called a resolution
sequence. Such sequences are constructed in
\cite{masur-minsky:complex2}, using what we call
``hierarchies of geodesics'' in $\CC(S)$. The machinery of 
\cite{masur-minsky:complex2} is cumbersome to describe fully, so we will 
give just a rough description of the construction and an indication of
how the results of \cite{masur-minsky:complex2} imply Theorem
\ref{Controlled Resolution Sequences}.

The construction is 
by an inductive procedure: we begin with what we call a ``tight
geodesic'', which is a sequence $m_i$ of simplices in $\CC(S)$ with
the first in $P$ and the last in $Q$, so that any sequence $\beta_i\in
m_i$ of vertices yields 
a geodesic in $\CC_1(S)$ joining $P$ to $Q$. (The $m_i$ satisfy an
additional condition which we need not use here).
The link of each $m_i$ is itself a
curve complex for a subsurface, in which $m_{i-1}$ and
 $m_{i+1}$ represent two simplices, and we construct a geodesic 
connecting them. We then repeat, with the complexity of the
subsurfaces decreasing at each step. (The actual construction 
is considerably
complicated by the need to take care of endpoints of geodesics
correctly, and by the fact that a typical simplex 
cuts up the surface into several components, in each of
which the construction continues independently).
The final structure is a collection of geodesics in subsurfaces
related by inclusion, and the pants decompositions $P_k$ are obtained
by taking ``slices'': picking a vertex at one level and then
inductively adding vertices from the geodesics supported in its
complementary subsurfaces. These slices can be ``resolved'' into the
sequence described in Theorem \ref{Controlled Resolution Sequences},
where the elementary moves $P_k\to P_{k+1}$
correspond to steps along geodesics in highest-level subsurfaces,
which are always one-holed tori or four-holed spheres. 

The fact that $J_\beta$ is always an interval follows from the proof
of Proposition 5.4 in \cite{masur-minsky:complex2}, in the course of
which we establish a monotonicity property of the way a resolution
steps through the geodesics in a hierarchy, that implies no curve is
ever repeated once it has been traversed. 

To obtain the inequality in part (3), note first that the length of
$J_{[s,t]}$ is just the sum of the lengths of the 
geodesics in the highest level subsurfaces meeting the slices 
based at $\beta_s,\ldots,\beta_t$. Each of these subsurfaces arises
from vertices in a lower-level (higher complexity) subsurface, so
their number is bounded by the sum of the lengths of the geodesics at
the lower level. Continuing inductively, if we have a length bound of
$B$ on all the geodesics encountered (and a length $(t-s)$ at the
bottom level), we obtain a bound of the form
$K(t-s)B^a$, where $a$ bounds the number of levels, which only depends
on the topological complexity of $S$. 

Finally, the length of a geodesic supported in a subsurface $Y$ is
bounded by a multiple of the projection distance $d_Y(P,Q)$: this is
the substance of Lemma 6.2 of 
\cite{masur-minsky:complex2}, which involves a crucial use of the
hyperbolicity property of $\CC(S)$, from \cite{masur-minsky:complex1}.
The bound of part (3) follows. 

Part (4) follows from the same construction, which in fact includes
annuli and their arc complexes as part of the discussion. The bound is
just a restatement of Lemma 6.2 of
\cite{masur-minsky:complex2} applied to annuli.

\section{The bounded shear lemma}
\label{shear}

In this section we will develop some estimates of shearing in annuli
that will be used near the end of the proof of the main theorem, in
Section \ref{final}.

We begin with an observation. 
Let $\sigma$ be a hyperbolic metric on
$S$, $\gamma$ a simple geodesic in $S$,  and
$B=\collar(\gamma,\sigma)$. Let $\hat Y$ be the compactified 
annular lift associated to $\gamma$ and $\hat B$ the annular lift of
$B$ to $\hat Y$.  Let $E$ denote one of the components of $\hat Y
\setminus \hat B$. 

There is only a bounded amount of twisting that a geodesic crossing
$\hat Y$ can do in $E$. In fact, let $\beta_1$ and $\beta_2$ be any
two geodesic lines connecting the two boundaries of $\hat Y$. 
We have:
\begin{equation}
  \label{bounded twist}
  d_{\AAA(E)}(\beta_1\intersect E,\beta_2\intersect E) \le 4.
\end{equation}

\realfig{bounded twist in B}{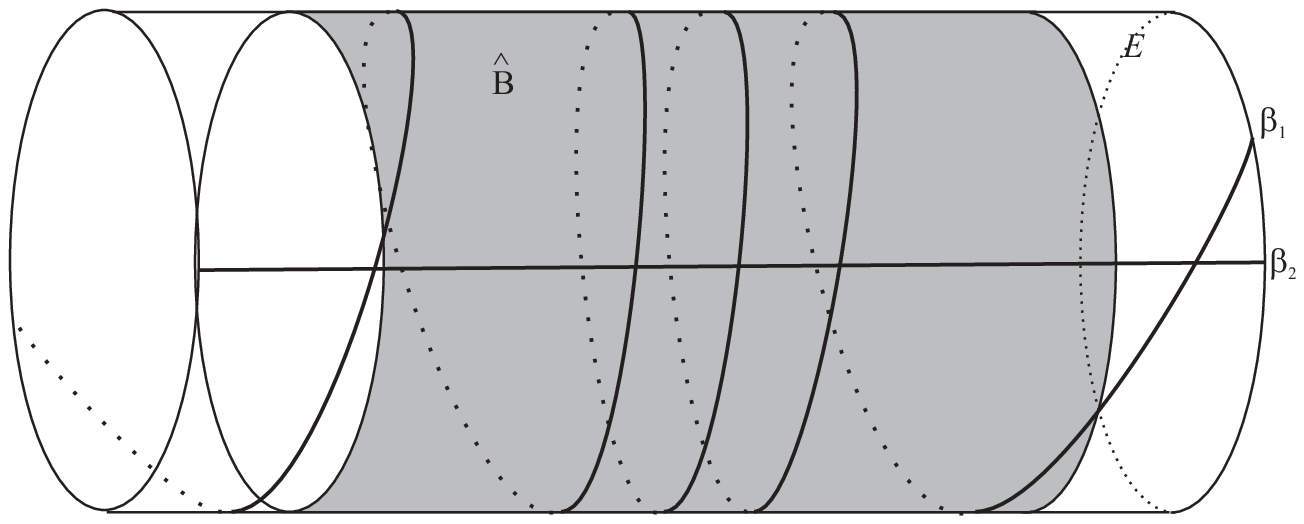}{The relative twisting
  $d_{\AAA(E)}$ of $\beta_1\intersect E$ and $\beta_2\intersect E$ is bounded.}

\begin{proof}
  The width  of the collar of $\gamma$ is the function
  $w(\ell(\gamma))$ described in the Appendix. 

  Lift $\hat B$ to a $w$-neighborhood $\til B$ of a geodesic lift
  $\Gamma$ of $\gamma$ in $\Hyp^2$. For $i=1,2$, $\beta_i$ lifts to a geodesic
  arc $\til\beta_i$ connecting $\boundary \til B$ to the circle at
  infinity. Let $p$ denote the length of its orthogonal projection to
  $\Gamma$. This is largest when $\til\beta_i$ is tangent to $\boundary
  \til B$, so let us consider this case. The arcs $\til\beta_i$ and its
  projection to $\Gamma$ form opposite sides of a quadrilateral with
  three right angles and an ideal vertex, whose two finite sides have
  lengths $w$ and $p$ (see Figure \ref{collar quad}). 
  Hyperbolic trigonometry (Buser \cite[2.3.1(i)]{buser:surfaces})
  gives us
$$
\sinh(w)\sinh(p) = 1.
$$
 On the other hand we have from the definitions in \S\ref{appendix}
 that $w=\max(w_0/2,w_0-1)$ where $w_0$ satisfies
$$
\sinh(w_0)\sinh(\ell/2) = 1.
$$
From this we obtain an expression for $p/\ell$ as a function of $w_0$,
and one can deduce 
\marginpar{VIA MAPLE!!}
that this quantity is bounded. Indeed the maximum
is obtained at $w_0=2$, and we have
$$
\frac{p}{\ell} < 1.5
$$ 
In other words, $\beta_i\intersect E$ travels less than 1.5 times around the
annulus, as measured by its orthogonal projection. 
We deduce that $\beta_1$ and $\beta_2$ intersect at most 3
times in $E$. The bound (\ref{bounded twist}) follows.
\end{proof}

\realfig{collar quad}{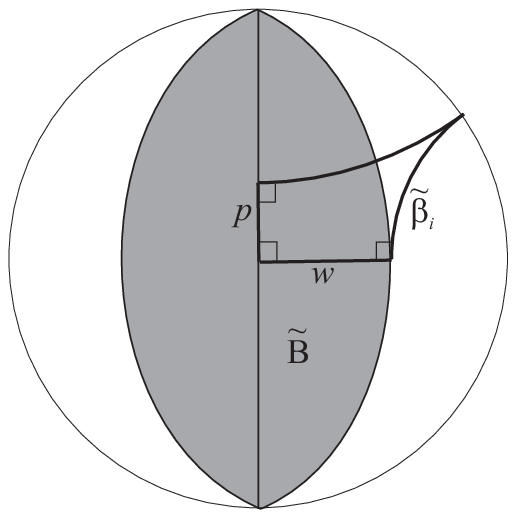}{The quadrilateral formed by
$\til\beta_i$ 
and its projection to the lift of $\gamma$}

\medskip

This observation prompts us to define the following measure of
shearing of two different metrics outside a collar. Let $\sigma,\tau$
be two hyperbolic metrics on $S$ and let $B$ be an annulus which is
equal to both $\collar(\gamma,\sigma)$ and
$\collar(\gamma,\tau)$. Lift to $\hat B$ in $\hat Y$ as above.
We let the {\em shear outside $B$} of the two metrics be the quantity
$$
\sup_{E,\alpha_\sigma,\alpha_\tau} d_{\AAA(E)}(\alpha_\sigma\intersect
E,\alpha_\tau\intersect E).
$$
Here $E$ varies over the two 
complementary annuli of $\hat B$ in $\hat Y$, $\alpha_\sigma$ varies
over all $\sigma$-geodesics that connect the two boundaries, and
$\alpha_\tau$ varies over all $\tau$-geodesics that connect the two
boundaries. 
Note that the
shear depends on the {\em pointwise} metrics, and not just their
isotopy classes.

The point of having a bound on this shear is that it allows us to
measure twisting by restricting to a collar. That is, suppose in the
setting above the shear outside $B$ is bounded by $D$. Then if
$\alpha_\sigma$ and $\alpha_\tau$ are any two $\sigma$ and $\tau$
geodesics, respectively, which cross $\gamma$, we immediately have
\begin{equation}\label{twists in collar}
|d_{\AAA(B)}(\alpha_\sigma\intersect B, \alpha_\tau\intersect B) 
- d_{\gamma}(\alpha_\sigma,\alpha_\tau)| 
\le 2D.
\end{equation}
Note that $d_{\AAA(B)}$ is a measure of twisting inside the collar
which depends on the particular curves we chose, whereas
$d_\gamma$ depends only on homotopy classes. 

The main lemma of this section bounds the shear outside a collar 
for pairs of metrics satisfying a special condition.

\begin{lemma+}{Shear bound}
Suppose $R$ is a subsurface of $S$ which is convex in two hyperbolic metrics
$\sigma$ and $\tau$, and that $\sigma$ and $\tau$ are locally
$K$-bilipschitz in the complement of $R$. Suppose that one component
of $R$ is an annulus $B$ which is equal to both
$\collar(\gamma,\sigma)$ and $\collar(\gamma,\tau)$ for a certain
curve $\gamma$. 

\realfig{lift sigma tau}{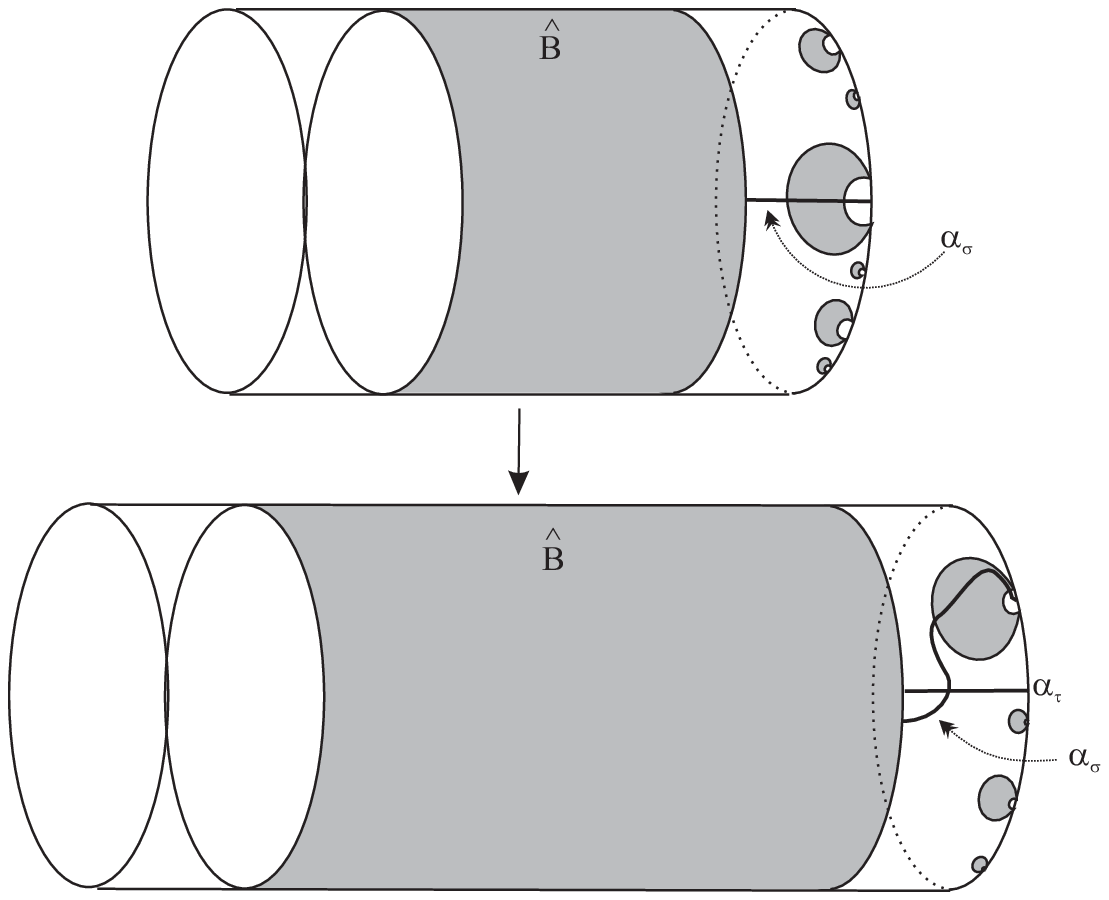}{
$\hat Y$, seen in the metric $\hat\sigma$ on top, and the metric $\hat\tau$
on the bottom. Some components of $\hat R$ are shown, shaded.}

Then the shear outside $B$ of $\sigma$ and $\tau$ is bounded above by
$\delta_0K$, where $\delta_0$ depends only on the topological type of $S$.
\end{lemma+}
\begin{proof}
Let $\hat Y$ be the compactified annular lift, and let $\hat\sigma$
and $\hat\tau$ be the lifts of $\sigma$ and $\tau$ to $\hat Y$.
Let $\hat R$ denote the lift of $R$ to $\hat Y$, and let $\hat B$ be
the component of $\hat R$ which is a homeomorphic lift of $B$.
Let $E$ be one of the
components of $\hat Y \setminus \hat B$. There
is a uniform upper bound $b$ to the $\hat \sigma$-length of the shortest
arc $\alpha$ connecting $\hat B$ to any component $D$ of $\hat R$
contained in $E$.  This comes from a standard area bound: 
Since $\ell_\sigma(\boundary B) $ is bounded below by $b_0$ (see
\S\ref{appendix}), 
an embedded collar of radius $r$ around $B$ has area at least
$b_0r$. Since the area of $S$ is fixed by the Gauss-Bonnet theorem,
$r$ is uniformly bounded. The first self-tangency of this collar
yields a bound for $b$.

Since by (\ref{bounded twist}) the choice of $\alpha_\sigma$
cannot change the twisting number in the lemma by more than 4, we
may assume that $\alpha_\sigma\intersect E$ has an initial segment
which is this 
shortest arc $\alpha$. Similarly we may assume that $\alpha_\tau$ is
orthogonal to $\boundary B$ in the metric $\tau$.
The Lipschitz condition in the
complement of $R$ implies that the $\hat\tau$-length of $\alpha$ is
bounded by 
$Kb$. Thus the number of essential intersections that
$\alpha_\tau$ can have with $\alpha_\sigma$ in the interval
$\alpha$ is bounded by $Kb/b_0$, since between any two such there is a
segment of $\alpha$ that projects orthogonally (in $\hat\tau$) to an
entire boundary component of $B$.
The remainder of $\alpha_\sigma$ is
contained in a convex set whose boundary $\xi$ is the lift of a boundary
component of $R$. Since $R$ is a convex subsurface
both in $\sigma$ and $\tau$,
$\xi$ still bounds a convex set in the metric $\hat \tau$. Thus
$\alpha_\tau$ can 
have at most one essential intersection with $\alpha_\sigma$ in this convex
set.  
%
%
\end{proof}

\section{The proof of the main theorem}
\label{final}

As discussed in the outline, after \cite{minsky:kgcc} it suffices to
prove this direction of  the theorem:
$$
\sup_Y d_Y(\nu_+,\nu_-) < K \quad\implies\quad
\inf_\gamma\ell_\rho(\gamma) > \ep
$$
for $\ep$ depending only on $K$ and $S$, with the infimum over
$\gamma$ that are not externally short.

\medskip

For simplicity of exposition, let us first prove the theorem in the
case when both ends of $N_\rho$ are degenerate. At the end of the
section we will indicate the changes in the argument necessary if one
or both of the ends are geometrically finite.

Let $K_1$ be a constant to be determined shortly, and
let $\ep_2$ be such that a $K_1$-neighborhood of any
$\ep_2$-Margulis tube is still contained in an $\ep_1$-Margulis tube
(see \S\ref{margulis tubes}).

Fix a closed curve $\gamma$ in $S$, and assume $\ell_\rho(\gamma)<
\ep_2$. In particular $\gamma$ must be simple (\S\ref{margulis tubes}),
and represents a vertex 
of $\CC_0(S)$. Our goal will be to bound the radius of the  Margulis
tube $\T_\gamma(\ep_0)$.

\segment{Initial pants}

Our first step is to obtain two pants decompositions
$P_+$ and $P_-$ of moderate $\rho$-length, and pleated surfaces
$f_\pm\in\pleat_\rho(P_\pm)$ which homologically encase $\T_\gamma(\ep_1)$.

\begin{lemma}{Encase gamma}
Suppose that $\rho\in\Dnp(S)$ has two degnerate ends, so that
$\nu_\pm\in\EL(S)$.

Let $\UU_+$ and $\UU_-$ denote neighborhoods
in $\UML(S)$ of $\nu_+$ and $\nu_-$,
respectively. There exist pants decompositions $P_+$ and
$P_-$ lying in $\UU_+$ and $\UU_-$, respectively, with the following
properties:
\begin{enumerate}
\item $\ell_\rho(P_\pm) \le L_1$
\item Given $f_+\in\pleat_\rho(P_+)$ and $f_-\in\pleat_\rho(P_-)$, 
$\T_\gamma(\ep_1)$ is homologically encased by $f_+$ and $f_-$.
\end{enumerate}
\end{lemma}

The last statement in the lemma means that $\T_\gamma(\ep_1)$ is covered with
degree 1 by any 3-chain whose boundary relative the cusps is $f_+ - f_-$.
In particular, any proper homotopy from $f_+$ to $f_-$ must cover all
of $\T_\gamma(\ep_1)$.

\begin{proof}
Since the $e_+$ end is degenerate, there is a sequence of pleated
surfaces $f^+_i$ eventually contained in every neighborhood of the end.
Choose $P^+_i\in\short(\sigma_{f^+_i})$. Then the geodesic
representatives of $P^+_i$ must also eventually exit $e_+$. 
\marginpar{REF -- I
have an argument somewhere} 
In particular they must converge to
$\nu_+$ in $\UML$. The same discussion applies with $+$ replaced by $-$.

Since $(N_\rho,cusps)$ is homeomorphic to $(S\times\R,cusps\times\R)$
and $\T_\gamma(\ep_1)$ is 
compact, any pleated surface in the proper homotopy class of $\rho$ which is
sufficiently far in a neighborhood of 
$e_+$ (or $e_-$) can be deformed to infinity, through cusp-preserving maps,
without meeting $\T_\gamma(\ep_1)$. 

Choosing $i$ high enough, then, insures that $f^+_i$ and $f^-_i$
homologically encase $\T_\gamma(\ep_1)$ and that $P^\pm_i\in\UU_\pm$.
Furthermore $P^\pm_i\in\short(\sigma_{f^+_i})$ implies that
$f^\pm_i\in\good_\rho(P^\pm_i, L_1)$, so that by Lemma \ref{Homotopy
  bound} there is a $K$-good homotopy (with $K$ depending on $L_1$)
between $f^\pm_i$ and any map in $\pleat(P^\pm_i)$. This homotopy has
bounded tracks outside the Margulis tubes of $P^\pm_i$, if any. Thus we may
choose $i$ high enough that this homotopy also avoids $\T_\gamma(\ep_1)$, 
so that $P_\pm = P^\pm_i$ satisfies the conclusions of the lemma.
\end{proof}

\segment{Resolution sequence and block map}

Join $P_+$ to $P_-$ with a resolution sequence 
$P_-=P_0\to\cdots\to P_n=P_+$, as in 
Theorem \ref{Controlled Resolution Sequences}.  Let
$\{\beta_i\}_{i=0}^m$ be the vertex sequence of the associated
geodesic in $\CC(S)$.  

Define a map $H:S\times[0,n]\to N_\rho$ as follows. For each
$j=0,\ldots,n$
choose $g_j\in\pleat_\rho(P_j)$, and for $j=0,\ldots,n-1$ let
$g_{j+\frac12}$ be a 
map in $\good_\rho(P_j,C_1)\intersect
\good_\rho(P_{j+1},C_1)$, promised to exist by Lemma \lref{Halfway
surfaces}. 

Lemma \lref{Homotopy bound} gives us a constant $K_1$ so that $g_j$
and $g_{j+\frac12}$ admit a $K_1$-good homotopy, and so do $g_{j+\frac12}$ and
$g_{j+1}$. Let $H|_{S\times[j,j+\frac12]}$ and
$H|_{S\times[j+\frac12,j+1]}$ be
these two
homotopies. Up to suitable precomposition by homeomorphisms of $S$
isotopic to the identity, we can arrange for the definitions to agree
on each integer and half-integer level, so they may be concatenated to
yield a single map $H$. The constant $K_1$ 
described here gives us our choice of $\ep_2$, which we note  depends
only on the topological type of $S$.

Let $\sigma(u) =\sigma_{H_u}$ be the induced metric on $S\times\{u\}$,
which we note is hyperbolic for $u$ an integer or half-integer. We
endow $S\times[0,n]$ with the metric that restricts to $\sigma(u)$ in
the horizontal direction and to the induced metric from $\R$ in the
vertical direction, and such that the  two directions are orthogonal.

\segment{Restricting the block map}

Assuming the neighborhoods $\UU_\pm$ have been chosen sufficiently small,
we can throw away all but a bounded number of the blocks and still
have the encasing condition. To see this, begin with the following claim:

\begin{claim}{intersection implies close}
There is a constant $C_2$ depending only on the topological type of $S$,
and a subinterval  $I_\gamma\subseteq[0,m]$ of diameter at most $C_2$, so
that $H(S\times[j,j+1])$ can meet $\T_\gamma(\ep_2)$ only if 
$j$ or $j+1$ are in $J_{I_\gamma}$.
\end{claim}

(The notation $J_I$ is defined in \S\ref{resolutions}.)

\begin{proof}
Let $\beta_i$ be a vertex that is in $P_j$.
If $H(S\times\{j\})$ meets $\T_\gamma(\ep_1)$ then (see
\S\ref{margulis tubes}) $\ell_{\sigma(j)}(\gamma) \le \ep_0$, 
and in particular  
$$\gamma\in\short(\sigma(j)) \subset
\Pi_\rho(\beta_i).$$ 
It follows
from Lemma \lref{Quasiconvexity} that 
$$
d_\CC(\beta_i,\gamma) \le C 
$$
where $C$ depends only on the topological type of $S$.
Thus, since $\{\beta_0,\ldots,\beta_m\}$ are the vertices of a geodesic,
the possible values of $i$ lie in an interval of  diameter at
most $2C$, which we call $I_\gamma$. In other words, $j\in J_{I_\gamma}$.

It remains to notice that, by the choice of $\ep_2$ and the $K_1$-good
homotopy property, if any part of a block $H(S\times[j,j+1])$ meets
$\T_\gamma(\ep_2)$ then one of the boundaries must meet
$\T_\gamma(\ep_1)$, and hence $j$ or $j+1$ are in $J_{I_\gamma}$.
\end{proof}

Let us therefore restrict our elementary
move sequence to 
$$P_{s-1}\to\cdots\to P_{t+1}$$
where $[s,t] = J_{I_\gamma}$.
This subsequence must still encase 
$\T_\gamma(\ep_2)$, since we have thrown away only blocks that avoid 
$\T_\gamma(\ep_2)$.
Let $M=t-s=|J_{I_\gamma}|$.

\medskip

In order to deduce a bound on $M$ from the diameter bound on
$I_\gamma$, we  consider
part (3) of  Theorem \ref{Controlled Resolution Sequences}, which tells us
that 
\begin{equation}\label{M bound}
        M \le K(2C)\sup_Y d_Y(P_+,P_-)^a,
\end{equation}
where the supremum is over subsurfaces $Y$ whose boundaries appear
among the $P_j$ in our subsequence. 
Such $P_j$ must lie in a $C+1$ neighborhood of $\gamma$, in the
$d_\CC$ metric. In order to compare $d_Y(P_+,P_-)$ to
$d_Y(\nu_+,\nu_-)$ for such $Y$, we will need the following lemma:

\begin{lemma}{projections converge}
There exists a constant $b$ depending only on $S$, so that
given $\gamma\in\CC_0(S)$, $R>0$ and $\nu\in\EL(S)$, there is a
neighborhood $\UU$ of $\nu$ in $\UML(S)$ for which the following holds: 

If $Y\subset S$ is a subsurface with $d_\CC(\gamma,\boundary Y) \le R$
and $\beta\in\CC_0(S)\intersect\UU$, then
$$
d_Y(\beta,\nu) \le b.
$$
\end{lemma}

\begin{proof}
By Klarreich's Theorem \ref{EL is boundary}, $\EL(S)$  is naturally
identified with the Gromov boundary of $\CC(S)$.
By the definition of the Gromov boundary, 
there is a neighborhood $\UU$ of $\nu$ such that, if
$\beta,\beta'\in\CC_0(S)$ are in $\UU$ then any geodesic $g$ in
$\CC_1(S)$ connecting them must lie outside an $R+2$-ball of $\gamma$. 
In particular every vertex of $g$ must be at least distance 2 from
$\boundary Y$, and hence must intersect $\boundary Y$ essentially.
Theorem 3.1 of \cite{masur-minsky:complex2}
states that in such a situation
$$
\diam_Y(\pi_Y(g)) \le A
$$ 
for a constant $A$ depending only on the topological type of $S$, and
in particular $d_Y(\beta,\beta') \le A$. 

Now consider a sequence $\beta_i$ converging to $\nu$ in $\UML(S)$,
with $\beta_0=\beta$. 
After restricting to a subsequence if necessary, the $\beta_i$ converge
in the Hausdorff topology 
to a lamination that includes the support of $\nu$,
which means that eventually $\diam_Y(\pi_Y(\beta_i)\union\pi_Y(\nu))
\le 2$. Thus  $d_Y(\beta,\nu)\le A+2$.
\end{proof}

Thus if we choose the original neighborhoods $\UU_+$ and $\UU_-$
sufficiently small, this lemma gives us
\begin{equation}\label{2b bound}
d_Y(P_+,P_-) \le d_Y(\nu_+,\nu_-) + 2b.
\end{equation}
The hypothesis of the main theorem bounds the right side, so together 
with (\ref{M bound}) we obtain
our desired uniform bound on $M$.

\medskip

Suppose that $\gamma$ is not a component of any $P_j$. 
Then according to Lemma \ref{Homotopy bound}, each block $H|_{S\times[j,j+1]}$
has track lengths of at most $2K_1$ within $\T_\gamma(\ep_1)$.
There are
only $M+2$ blocks  in our restricted sequence, and they cover all of
$\T_\gamma(\ep_1)$. The beginning and the end of the sequence are
outside of $\T_\gamma(\ep_1)$, so that any point  in 
$\T_\gamma(\ep_1)$ is at most
$$2K_1(M+2)$$ 
from its boundary. This bounds the radius of $\T_\gamma(\ep_1)$ by 
$2K_1(M+2)$, and a corresponding lower bound for
$\ell_\rho(\gamma)$  follows as in \S\ref{margulis tubes}.

\segment{Bounding twists}

Now suppose that $\gamma$ does appear among the $\{P_j\}$. Then
$J_\gamma$ is some nonempty subinterval of $J_{I_\gamma}$ by
Theorem \ref{Controlled Resolution Sequences}, and we let $\alpha$ and
$\alpha'$ be the predecessor and successor curves to $\gamma$ in the
sequence. Both of them cross $\gamma$, and we have
by part (4) of Theorem \ref{Controlled Resolution Sequences} that
$d_\gamma(\alpha,\alpha')$ is uniformly approximated by
$d_\gamma(P_+,P_-)$, which by (\ref{2b bound}) and the hypothesis of
the main theorem, is uniformly bounded. 
Let $D$ denote this bound.

Write $J_\gamma = [k,l]$.
By the normalization used in Lemma
\ref{Homotopy bound}, for all integer and half-integer $u\in [k,l]$
the annuli $\collar(\gamma,\sigma(u))$ coincide. Name this common
annulus $B$. 
consider the solid torus 
$$
U = B\times[k-\half,l+\half].
$$
The map $H$ can take the complement of $U$
at most $2K_1(M+2)$  into $\T_\gamma(\ep_1)$, by the same argument as above,
and hence there is a uniform
$\ep_3>0$ so that $H(U)$ must cover
$\T_\gamma(\ep_3)$. Assume $\ell_\rho(\gamma)<\ep_3$, for otherwise we
are done. 

Consider the geometry of $\boundary U$, in the metric we have placed
on $S\times[0,n]$. The top annulus
$B\times\{k-\half\}$ is $\collar(\gamma,\sigma({k-\half}))$.
Since $H_{k-\frac12}\in\good(\gamma,C_1)$ by construction, we have an
upper bound of $\ep_3+C_1$ on the circumference of this collar. 
We also have a lower bound of $\ep_3$ since the image of
$H_{k-\frac12}$ avoids the interior of $\T_\gamma(\ep_3)$.
Thus, this annulus has uniformly bounded geometry. 
The same holds for $B\times\{l+\half\}$ in the metric 
$\sigma({l+\half})$

The vertical annuli $\boundary B \times [k-\half,l+\half]$
have height bounded by
$M+1$. Thus $\boundary U$ is a concatenation of four annuli of bounded
geometry, and hence is itself
up to uniformly bounded distortion a Euclidean square torus, and $H$ maps
it with a uniform Lipschitz constant into $N$. 
Let us now try to control a meridian curve for $U$.

Realize $\alpha$ as a geodesic in $\sigma({k-\half})$
and let $a$ be an arc of $\alpha\intersect B$ (there may be two). Similarly
assume $\alpha'$ is a geodesic in $\sigma({l+\half})$ 
and let $a'$ be an arc of $\alpha'\intersect B$. 
Lemma \ref{halfway annulus bound} gives an upper bound for the length
of $a$ in $\sigma({k-\half})$, and for the length of 
$a'$ in $\sigma({l+\half})$.
The curve
$$
m = \boundary (a\times[k-\half,l+\half]),
$$
is a meridian of $U$, and we claim that its length is uniformly bounded.

\marginpar{NEED PICTURE}

The arc $a$ may a priori be
long in $\sigma({l+\half})$,
but its length is estimated by the number
of times it twists around the annulus, which in turn is estimated by
$d_{\AAA(B)}(a,a')$ since $a'$ has bounded length in this metric.

We can bound this twisting using the results in Section \ref{shear}.
Note first that, for each $j=k-\half,k,\ldots,l-\half,l$,
the pair of metrics $\sigma(j),\sigma({j+\half})$
satisfy the hypotheses of Lemma \lref{Shear bound}. That is, they are
$K_1$-bilipschitz outside a union of collars including $B$.
Thus, Lemma \ref{Shear bound} bounds the shearing
outside $B$ at each step, and the number of steps is bounded by
$2M+2$. We conclude that there is a bound 
on  the shearing outside $B$ of the two metrics
$\sigma(k-\half)$ and $\sigma(l+\half)$.
It follows as in (\ref{twists in collar}) that we have a bound of the form
\begin{equation}\label{measure twist in collar}
|d_{\AAA(B)}(a,a') - d_\gamma(\alpha,\alpha')| \le M'
\end{equation}

With this estimate and the bound  $d_\gamma(\alpha,\alpha')\le D$,
we find that $a$
and $a'$ intersect a bounded number of times, so that the length of
$a$ is uniformly bounded in $S\times\{l+\half\}$. It follows that 
$H(m)$ is uniformly bounded. 
It therefore spans a disk of bounded diameter,
and now by a coning argument
we can homotope $H$ on all of $U$ to a new map of bounded
diameter. This bounds the radius of $\T_\gamma(\ep_3)$ from above,
and again we are done.

\segment{The case of geometrically finite ends}

The main change in the argument comes at the beginning, in our choice
of $P_+$ and $P_-$. 


Suppose that $e_+$ is geometrically finite. It suffices to consider 
$\gamma$ which is not externally short, i.e. $\ell_{\nu_+}(\gamma) > \ep_0$,
because the second statement in the main theorem considers only such
curves, and the first statement is insensitive to the removal of a
finite number of curves from consideration. The choice of $\ep_1$ and
Sullivan's theorem (see \S\ref{margulis tubes}) imply that we also have 
$\ell_{\nu'_+}(\gamma) > \ep_1$.

We may therefore choose a pants
decomposition  $P_+\in\short(\nu'_+)$ so that $\gamma$ is not one of
its components (our choice of $L_1$ was made so that this would be
possible, see \S\ref{projection coefficients}).
Let $f'_+: S \to \boundary_+(C(N_\rho))$ be the pleated surface in the
homotopy class of $\rho$ that parametrizes the $e_+$-boundary  of the
convex hull. In particular $f'_+\in\good_\rho(P_+,L_1)$.
Thus, choosing $f_+\in\pleat_\rho(P_+)$, 
Lemma \ref{Homotopy bound} gives  us a $K_2$-good homotopy $G_+$ between
$f'_+$ and $f_+$, with $K_2$ depending on $L_1$. Since $\gamma$ is not
a component of $P_+$, $G_+$ can only penetrate a distance $K_2$ into
$\T_\gamma(\ep_1)$, and hence does not meet $\T_\gamma(\ep_4)$, 
with $\ep_4$ depending only on $K_2$.

We do the same thing with $f_-,P_-$ if $e_-$ is geometrically finite,
or we repeat the original discussion if $e_-$ is degenerate. Thus we
have $f_+,f_-$ encasing $\T_{\gamma}(\ep_4)$.

We can therefore continue as before, constructing a resolution
sequence from $P_-$ to $P_+$ and then restricting it. The argument
goes through with slightly altered constants (since we have replaced
$\ep_1$ with $\ep_4$), and the one step that needs attention is the
comparison (\ref{2b bound}) between $d_Y(P_+,P_-)$ and
$d_Y(\nu_+,\nu_-)$, for appropriate $Y$. In other words we must bound
$$
d_Y(P_+,\nu_+)$$
in the geometrically finite case (and similarly for $e_-$). 
This quantity by definition is just $\min_{Q_+}d_Y(P_+,Q_+)$ where
$Q_+$ varies over pants decompositions
in $\short(\nu_+)$, provided $P_+$ meets $Y$ nontrivially
and $Q_+$ can be found which does the same.

Since $\nu_+$ and $\nu'_+$ admit a $K_0$-bilipschitz map by Sullivan's
theorem (\S\ref{end invariants}),
the $\nu'_+$-length of $Q_+$ is bounded by $K_0L_1$, which gives a
bound on its 
intersection number with $P_+$. This bounds $d_Y(P_+,Q_+)$ 
provided both pants decompositions intersect $Y$.

If $Y$ is not an annulus (and recall that we never consider
three-holed spheres), then it automatically intersects any pants
decomposition.  Thus the only problematic case is if $Y$ is an annulus
whose core is a component of $P_+$ or $Q_+$.
However, we recall that the only annulus that actually comes into the
argument is the annulus with core $\gamma$. 
Since $\gamma$ is not externally short, we have already chosen $P_+$
so that $\gamma$ is not a component of it, and we can do the same for $Q_+$.

The rest of the proof goes through in the same way.
\qed
\section{Appendix: Hyperbolic geometry}
\label{appendix}

In this appendix we write out statements, and sketch some proofs, of a
few facts and constructions in hyperbolic geometry. 
These are ``well-known'' in the sense that those working in this field
are familiar at least with some variation of them, or would find it
straightforward to derive them. Still it seems advisable to include
some discussion.

Throughout, a {\em hyperbolic surface} always means a finite area
hyperbolic surface which could have closed geodesic boundary
components and/or cusps. By abuse of notation we usually think of a
cusp as a boundary component of length 0.

\segment{Collars}

Let $w_0$ be the function
$$
w_0(t) = \sinh^{-1}\left(\frac{1}{\sinh\left(t/2\right)}\right).
$$
For a simple closed geodesic $\gamma$ of length $\ell$ in a
hyperbolic surface 
$X$, we define 
$$
\collar_0(\gamma,X) =
\{p\in X: \dist(p,\gamma) \le w_0(\ell)\}.
$$
When the ambient surface $X$ is understood we omit it from the notation.
This set is always an embedded annulus, and in fact if
$\gamma_1\ldots,\gamma_k$ are disjoint and  homotopically distinct
then $\collar_0(\gamma_i)$ are pairwise disjoint. 
(See e.g. Buser \cite[Chapter 4]{buser:surfaces}.)

We will need a slightly smaller collar, so as to guarantee a definite
amount of space in its complement. Define
$$
w = \max(w_0/2,w_0-1)
$$
and let 
$$\collar(\gamma,X) = \{p\in X: \dist(p,\gamma) \le w(\ell)\}.
$$

Let $\gamma'_0$ be a boundary component of $\collar_0(\gamma)$
(assume $\gamma'_0\ne\gamma$ if $\gamma$ itself is in $\boundary X$),
and similarly let $\gamma'\ne \gamma$ be a boundary component of $\collar(\gamma)$.
The length of $\gamma'_0$ (resp. $\gamma'$)
is given by  $\ell\cosh(w_0(\ell))$ (resp. $\ell\cosh(w(\ell))$). A bit of
arithmetic shows that
\marginpar{Arithmetic $=$ Maple}
\begin{equation}\label{collar boundary bound}
l(\gamma') \le l(\gamma'_0) \le \ell(\gamma) + 2
\end{equation}

We can define $\collar_0(\gamma)$ and $\collar(\gamma)$ also when
$\gamma $ represents a cusp of $X$, 
describing them either explicitly or
as a limit as $\ell\to 0$. 
The boundary of $\collar_0(\gamma)$ for a cusp is horocyclic, and its length 
is 2 (the limiting value of $\ell\cosh(w_0(\ell))$ as $\ell\to 0$).
The slightly smaller $\collar(\gamma)$ has horocylic boundary a
distance 1 inside the boundary of $\collar_0(\gamma)$, so its length
is $2/e$.

We remark in fact that the boundary length $\ell\cosh w_0(\ell)$ of
$\collar_0(\gamma)$ 
is increasing with $\ell$, and hence always at least 2. There is a
similar lower bound $b_0$ for the boundary of $\collar(\gamma)$, which
we will not compute explicitly.

If $P=\{\gamma_1,\ldots,\gamma_k\}$ is a curve system we let
$\collar(P) = \union_i\collar(\gamma_i)$.

\segment{Hyperbolic pairs of pants}

If $Y$ is a hyperbolic surface as above, with genus 0 and
three boundary components, we call it a hyperbolic pair of pants. 
The three boundary lengths determine the metric on $Y$ completely (up to
isotopy). Let $Y_0 = Y\setminus int(\collar(\boundary Y))$.

A {\em tripod} is a copy of the 1-complex $\Delta$ obtained from three
disjoint copies of $[0,1]$ (called ``legs'') by identifying the three
copies of 
$\{0\}$. The three copies of $\{1\}$ are called the boundary of
$\Delta$. An {\em essential tripod} in $Y_0$
is an embedding of $\Delta$ (also called $\Delta$) taking
$\boundary \Delta$ to $\boundary Y_0$, such that each subarc of
$\Delta$ obtained by deleting one copy of $(0,1]$ is not homotopic rel
endpoints into $\boundary Y_0$ (Figure \ref{tripod examples}).

\realfig{tripod examples}{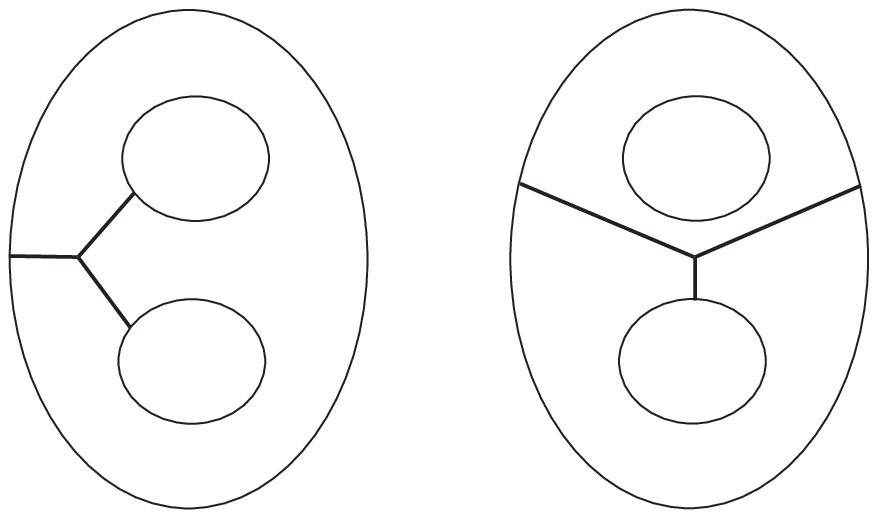}{The two types of essential tripods}

\begin{lemma}{Short tripods}
There is a constant $\delta>0$ such that for
any hyperbolic pair of pants $Y$ and each boundary component $\gamma$ of $Y_0$,
there is an essential tripod
$\Delta\subset Y_0$ whose three legs have length
at most $\delta$, and which meets $\gamma$.

On the other hand, there exists $\delta'>0$ such  no essential
tripod in any $Y_0$ has all three legs of length less than $\delta'$.
\end{lemma}

The next lemma describes how an estimate on the boundary lengths of
$Y$ yields an estimate on its isometry type.

\begin{lemma}{pants bilipschitz}
Given $C\ge 0$ there exists $K\ge 1$ such that the following holds. 

Let $Y$ and $Y'$ be two hyperbolic pairs of pants with boundaries 
$\{\gamma_i\}_{i=1}^3$ and $\{\gamma'_i\}_{i=1}^3$, respectively.
Suppose that 
$$
|\ell(\gamma_i)-\ell(\gamma'_i)| \le C.
$$
Then there is a homeomorphism $h:int(Y)\to int(Y')$
taking $Y_0$ to $Y'_0$ and 
$\collar(\gamma_i)$ to $\collar(\gamma'_i)$, which is $K$-bilipschitz
on $Y_0$.
Furthermore, for any arc $a$ on $\boundary Y_0$ we have
$|l(\alpha) - l(h(\alpha))| \le K$.
\end{lemma}

Let us prove Lemma \ref{pants bilipschitz} first. 
Any hyperbolic pair of pants $Y$ admits a canonical
decomposition into two congruent right-angled hexagons, by cutting
along the shortest arcs connecting each pair of boundary components
(see Buser \cite{buser:surfaces}. In the limiting case of cusps we obtain
degenerate hexagons with ideal vertices). Thus it will suffice to
find appropriate bilipschitz maps between such hexagons.

\newcommand\ra[2]{r_{#1#2}}
\newcommand\Ba[2]{B_{#1#2}}  
\newcommand\hex{\Xi}

\segment{Comparison of hexagons}

Let $\hex$ be a hyperbolic right-angled hexagon, with three
alternating sides $A_1,A_2,A_3$ of lengths $a_1,a_2,a_3$. 
(we will adopt the convention here of denoting edges by
capital letters and their lengths by the corresponding lower case
letters). 
Each side $A_i$ admits an embedded ``collar rectangle'', which we call
$\collar(A_i,\hex)$, and is just the $w(2a_i)$-neighborhood of $A_i$
in $\hex$. If $Y$ is the hyperbolic pair of pants obtained from $\hex$
by doubling along its other three boundaries, then $\collar(A_i,\hex)$
is clearly just $\collar(\gamma_i,Y)\intersect \hex$, where $\gamma_i$
is the double of $A_i$.

Lemma \ref{pants bilipschitz} clearly follows from the following
property of hexagons:

\begin{lemma}{hexagon comparison}
Given $C>0$ there exists $K\ge 1$ so that the following holds:

Let $\hex,\hex'$ be two hyperbolic right-angled hexagons with
alternating sides $\{A_i\}$ and $\{A'_i\}$, respectively.
Suppose that
$$
|a_i-a'_i| \le C, \quad i=1,2,3.
$$
Then there is a label-preserving homeomorphism
$f:\hex\to\hex'$ which takes $\collar(A_i,\hex)$ to $\collar(A'_i,\hex')$
and is $K$-bilipschitz on the complement of the collars.
Furthermore for any arc $\alpha$ in $\union A_i$, we have
$|l(\alpha) - l(f(\alpha))| \le K$.
\end{lemma}

\begin{proof}
We begin by describing a decomposition of $\hex$ into controlled
pieces in a way that is determined by the $\{a_i\}$. 
In the following,  $(i,j,k)$ always denotes a permutation of
$(1,2,3)$.

\subsection*{Case 1:} Suppose that the three ``triangle inequalities''
\begin{equation}
  \label{eq:triangle}
  a_i \le a_j + a_k
\end{equation}
all hold.  Then there is a unique triple $\ra12,\ra23,\ra13\ge 0$ satisfying
\begin{equation}
  \label{eq:r def}
  \ra ij + \ra ik = a_i.
\end{equation}
(where we use the convention $\ra ij = \ra ji$). In fact we simply set 
$\ra ij = \half(a_i + a_j - a_k)$.

Now define three ``bands'' $\Ba ij$ in $\hex$ as follows: 
$\Ba ij$ is the (closed) $\ra ij$-neighborhood of the edge $C_{ij}$ 
of $\hex$  which is the common perpendicular of $A_i$ and $A_j$. See Figure
\ref{Bands 1}.

\realfig{Bands 1}{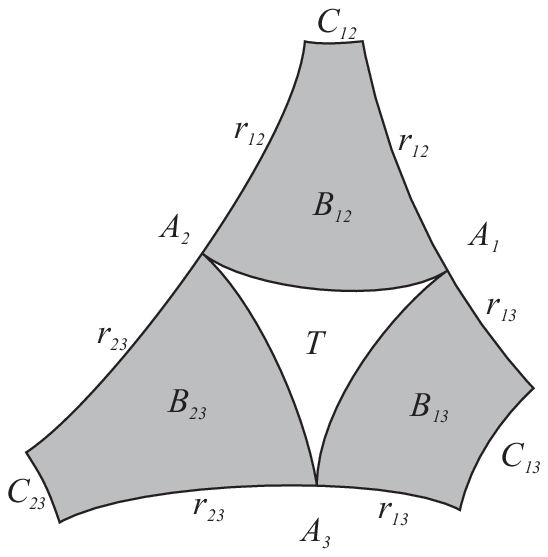}{The bands decomposition of $\hex$ in
  case 1.}

The two segments $\Ba ij\intersect A_i$ and $\Ba ik\intersect A_i$
cover all of $A_i$ and meet in their common boundary point. The
interiors of $B_{ij}$ are disjoint as a consequence of the triangle
inequality and the fact that $\dist(C_{ij},C_{ik}) = a_i$.
The closure $T$ of the complement $\hex\setminus \bigcup_{ij} \Ba ij$
is a ``triangle'' with the following properties: 

\begin{itemize}
\item
Each edge $E_{ij}= T\intersect \Ba ij$ has curvature $\kappa_{ij}\in[0,1]$,
which is in fact given by $\tanh \ra ij$. The curvature vector points outward.
\item
Adjacent edges meet at angle 0.
\end{itemize}
From this and the Gauss-Bonnet theorem we 
have 
$$ \sum_{ij=12,13,23} e_{ij}\kappa_{ij} = \pi - \mathrm{area}(T)$$
which implies the upper bound
\begin{equation}
\label{T edge upper}
e_{ij} \le \frac{\pi}{\tanh \ra ij}.
\end{equation}
There is also a lower bound
\begin{equation}
\label{T edge lower}
e_{ij} \ge 2\log \frac2{\sqrt3}
\end{equation}
obtained as follows: Extend each $E_{ij}$ to a line of constant
curvature in $\Hyp^2$. Each bounds a region $R_{ij}$ that contains
$B_{ij}$, so that the three regions have disjoint interiors. Since
$\kappa_{ij} \le 1$, through any point 
$x\in\boundary R_{ij}$ passes the boundary of a horoball
contained in $R_{ij}$. One can check that for any three horoballs in
$\Hyp^2$ with disjoint interiors, and any point $p$, the distance from
$p$ to at least one of the horoballs is at least $\log \frac2{\sqrt3}$ (the
extreme case is where all three are tangent). The midpoint of $E_{ij}$ 
meets a horoball in $R_{ij}$ and
is a distance at most $e_{ij}/2$ from horoballs in $R_{ij}$ and $R_{jk}$.
Thus we have (\ref{T edge lower}).

The geometry of the band $\Ba ij$ is easy to describe: it can
parametrized by the rectangle $[0,\ra ij]\times[0,c_{ij}]$ with the
metric
\begin{equation}\label{band 1 metric}
dx^2 + \cosh^2(x)dy^2,
\end{equation}
Where $C_{ij}$ is identified with $\{0\}\times [0,c_{ij}]$.

\subsection*{Case 2:} Suppose that one of the opposite triangle inequalities
holds, e.g.
\begin{equation}
  \label{eq:untriangle}
  a_1 \ge a_2 + a_3.
\end{equation}
 
We then let $\ra 12=a_2$, $\ra 13 = a_3$, and 
$$
\ra 11 = a_1-a_2-a_3,
$$
So that $a_1 = \ra 12 + \ra 11 + \ra 13$. We then have bands 
$\Ba 12$ and $\Ba 13$ as before, and
$\Ba 11$ is defined as follows (see Figure \ref{Bands 2}):
Let $H_1$ be the common perpendicular of $A_1$ and its opposite
edge $C_{23}$. In $A_1$ let $J_{11}$ be the closure of the complement of 
$\Ba 12 \intersect A_1$ and
$\Ba 13 \intersect A_1$. This has length $\ra 11$. Join each $x\in
J_{11}$ to a point $y\in C_{23}$ with a curve equidistant from $H_1$.
We obtain a foliated rectangle, $\Ba11$ (if $\ra11=0$ then $\Ba11$ is
a single segment). Similarly to the other bands, we can describe the
metric in $\Ba11$ by the formula (\ref{band 1 metric}), but on a
rectangle of the form
$ [u_2,u_3]\times[0,h_1]$, where $u_3-u_2 = \ra 11$.

\realfig{Bands 2}{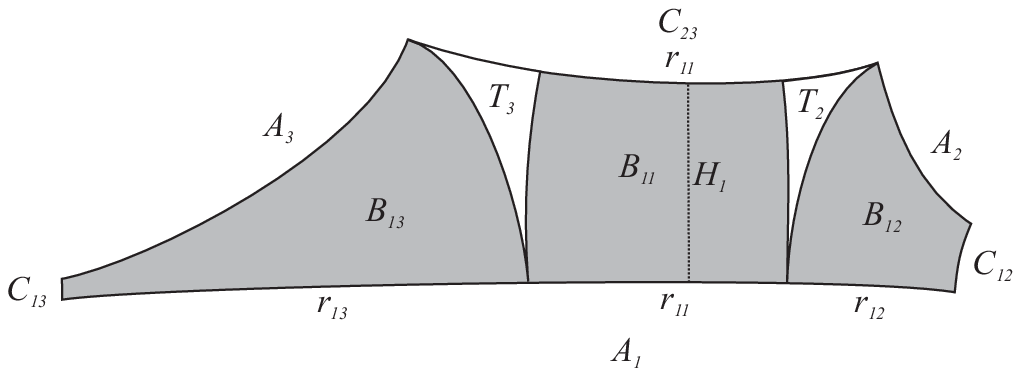}{The bands decomposition of $\hex$ in
  case 2.}

The closure of $\hex \setminus \bigcup_i \Ba1i$ is now two triangles
$T_2,T_3$, with
angles $0$, $0$ and $\pi/2$. Let $E_{1m}=T_m\intersect B_{1m}$ and let
$E_{11,m} = T_m\intersect B_{11}$, for $m=2,3$. Note that the
curvature of $E_{11,m}$ is $\kappa_{11,m} = \tanh|u_m|$.
Upper and lower bounds on $e_{1m}$ and $e_{11,m}$ can be derived from 
(\ref{T edge upper}) and 
(\ref{T edge lower}), by observing that the union of $T_m$ with its
reflection across $C_{23}$ is a triangle of the same type as $T$ in
case 1. For $e_{1m}$ we obtain (\ref{T edge upper}) and 
(\ref{T edge lower}) exactly, whereas for $e_{11,m}$ we get
\begin{equation}\label{T edge upper 2}
e_{11,m} \le \frac{\pi}{2\tanh|u_m|}
\end{equation}
and
\begin{equation}
  \label{T edge lower 2}
  e_{11,m} \ge \log \frac2{\sqrt3}.
\end{equation}
Let $|u_m|$ be the larger of $|u_2|,|u_3|$. 
Then $|u_m|\ge r_{11}/2$, and since $e_{11,m}$ must then be the larger of
$|e_{11,2}|,|e_{11,3}|, $ we have a bound for $j=2,3$:
\begin{equation}\label{T edge upper 2 both}
e_{11,j} \le \frac{\pi}{2\tanh(r_{11}/2)}
\end{equation}

\subsection*{The comparison:}
Now consider two hexagons $\hex, \hex'$ satisfying the 
bound $|a_i-a'_i|\le C$. 

Suppose first that both $\hex$ and $\hex'$ are in case 1. 
Note immediately that we also have $|\ra ij - \ra ij'| \le
\frac{3C}{2}$.

If $r_{ij} > 3C$ for each $i,j$, then $r'_{ij} > \frac{3C}{2}$.
In this case we have upper bounds of $\pi/\tanh(3C/2)$ on both
$e_{ij}$ and $e'_{ij}$. The bilipschitz map $\hex\to\hex'$ can be 
constructed separately on each piece of the decomposition.

Map each band $B_{ij}$ to $B'_{ij}$ using an affine stretch 
on the parameter rectangles.
The metric described in (\ref{band 1 metric}) gives 
a uniform bilipschitz bound on this map.
Note that the assumption that both $\ra
ij$ and $r'_{ij}$ are bounded away from 0 bounds the bilipschitz
constant in the $x$ direction, and the fact that their difference is
bounded gives a bound in the other direction (since the $\cosh$ factor
gives an exponential scaling).

The triangles $T$ and $T'$ are also in uniform bilipschitz
correspondence, since in fact they vary in a compact family of
possible figures (due to the length and curvature bounds). 

\medskip

If at least one $r_{ij}<C$ we must take a bit more care. Let us
consider a limiting case where at least one $r_{ij} = 0$ and the rest
are no smaller than $C$. These cases are illustrated in Figure
\ref{case 1 boundary}. 

\realfig{case 1 boundary}{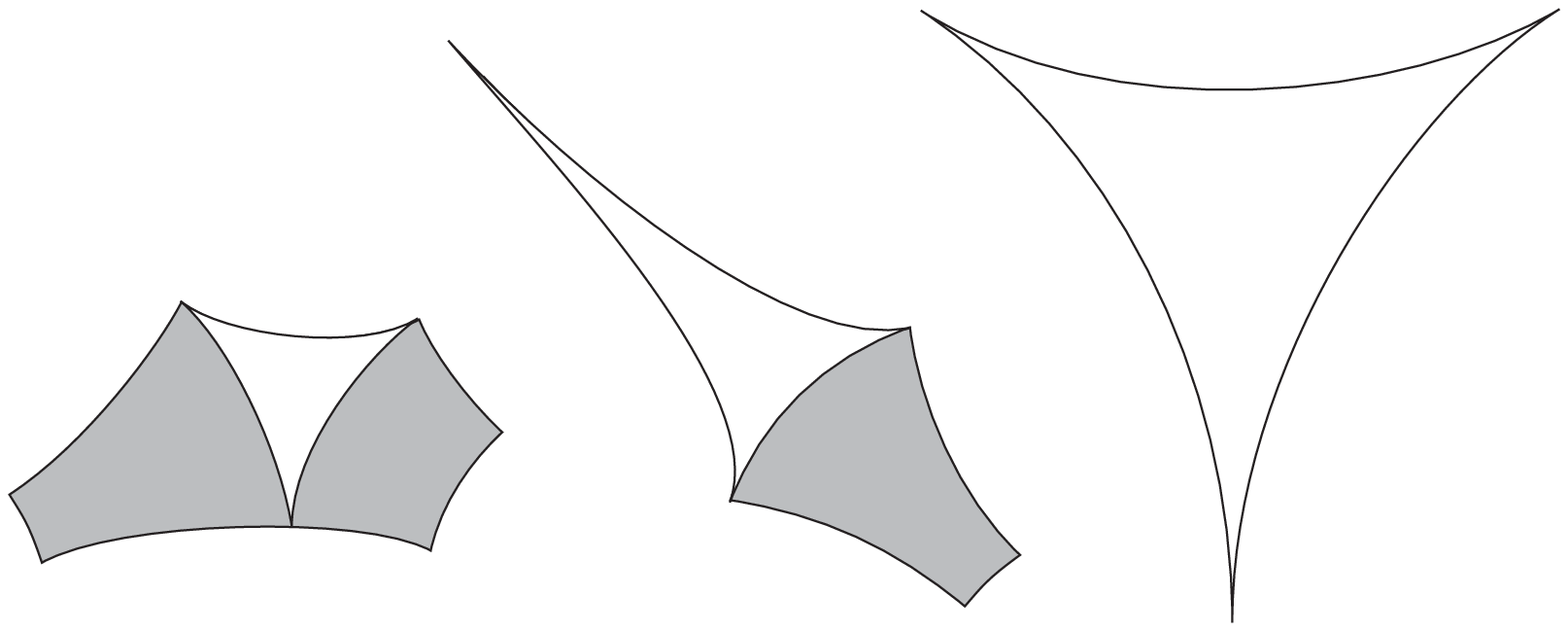}{The boundary cases where some
$r_{ij}=0$. From left to right, one, two and all three are 0.}

If $\ra12=\ra23=\ra13=0$ then $a_1=a_2=a_3=0$ and $\hex$ is an ideal
triangle. Recall that we are interested only in bounds on the
complements of the collars, and this is in this case a right-angled ``hexagon''
with three of the sides equal to horocyclic arcs of definite length.
If $\ra12=\ra13=0$ but $\ra23>C$, then $T$ has two geodesic edges
meeting in an ideal vertex, and a third leg $E_{23}$ satsifying the
bounds (\ref{T edge upper},\ref{T edge lower}). To this leg is
attached the band $B_{23}$.

If $\ra12=0$ and $\ra13,\ra23>C$ then we obtain a triangle $T$ with one
geodesic leg and two curved legs to which are attached the bands. Note
that the length bounds on the curved legs imply a length bound on the
straight leg.

If we now take a general $\hex$ with some $r_{ij}<C$, consider the
family of shapes obtained by taking those $r_{ij}$ to 0, and fixing
the rest. These must converge to one of the three cases above, and
the bands corresponding to $r_{jk}\ge C$ remain fixed. Thus we can
find some uniformly bilipschitz mapping taking $\hex$ to the limiting
case. \marginpar{Note it will not exactly preserve the long bands -- discuss?}

Two hexagons $\hex$ and $\hex'$ both near one of the boundary cases
are therefore near to each other, via the composition of the two maps.

\medskip

Now suppose our two hexagons are in case 2. Again we obtain bounds
$|r_{1j}-r'_{1j}|\le C$ for $j=2,3$ and
$|r_{11}-r'_{11}| \le 3C$. Assume first $r_{1j}>6C$ for $j=1,2,3$.
The bands $B_{12},B_{13}$
again have uniformly bilipschitz maps to $B'_{12}$ and $B'_{13}$,
respectively, coming from the affine stretches of the parameter
rectangles. The same is true for $B_{11}\to B'_{11}$, using the 
parametrization $[u_2,u_3]\times[0,h_1]$ with metric (\ref{band 1
metric}), together with the bounds (\ref{T edge upper 2}--\ref{T edge
upper 2 both}).

The cases where some $r_{1j}\le 6C$ are handled as before, by
comparing to the boundary cases where some $r_{1j}=0$. 
We leave the details to the reader.

\medskip

It is also possible for $\hex$ to be in  case 1 and $\hex'$ to be in
case 2. However, this can only happen if both of them are near their
respective boundary cases, and one can check that these boundary cases
are in fact in the overlap of case 1 and case 2, and hence can be
compared with each other.

\medskip

What we have shown is that $\hex $ and $\hex'$ admit a homeomorphism
which is bilipschitz outside the collars of those $A_i,A'_i$ that are
sufficiently short. For the long $A_i$ we observe that, 
as in (\ref{collar boundary bound}),  the length of the boundary of 
$\collar(A_i)$ in $int(\hex)$ is at most
$a_i+1$, and the collar width $w$ has been chosen
so that it is at most half the width of the largest embedded
collar. It follows from this that the map can be adjusted to take all
collars of $\hex$ to collars of $\hex'$, with a bounded change in the
bilipschitz constant. The additive distortion in length of all subarcs
of $A_i$ follows from the construction.
\end{proof}

\medskip

The proof of Lemma \ref{Short tripods} is fairly simple from the
geometric picture of hexagons we have just produced. 
Let $Y$ be decomposed as the double of a hexagon $\hex$. If $\hex$ is
in Case 1, the curved triangle $T$, minus the collars of the $A_i$, is
a bounded-diameter hexagon, and an essential tripod of bounded size
can be embedded in it. This tripod meets all three
boundary collars.

In Case 2, we consider the double of $T_2$ (or $T_3$) across
$C_{23}$. This is a curved triangle of the same type as $T$, and after
cutting away the intersections with the collars we can again obtain a
bounded essential tripod. Note that this tripod has two legs
terminating on the same boundary component of $Y$. It only meets two
boundary collars, but since we can  choose both $T_2$ and
$T_3$, for each collar there is a tripod meeting it.

To obtain the lower bound on any essential tripod in $Y$, we argue as
in the lower bound on $e_{ij}$ in the previous proof. Suppose the legs
of the tripod have length less than $l$. 
Lifting the tripod to $\Hyp^2$ we see that the intersection point of
the legs has  distance $\le l$ from three distinct lifts of the collars,
each of which is bounded by by a curve of constant curvature
$\kappa\in[0,1]$. The smallest possible $l$ is again obtained if
$\kappa=1$ for each curve, namely $l\ge 2\log\frac2{\sqrt3}$.
\qed

\segment{A curve shortening lemma}

Let $\gamma$ be a circle of length $l(\gamma)$,
and let $h:\gamma\to N$ be a homotopically nontrivial,  arclength-preserving
map into a hyperbolic 3-manifold $N$. Lift $h$ to a map  $\til
h:\til\gamma\to\Hyp^3$. After appropriate choice of basepoints, 
$h_*(\pi_1(\gamma))$ is a nontrivial cyclic subgroup of $\pi_1(N)$ preserving
$\til h(\til\gamma)$. Let $\ell_N(h)\ge 0$ denote the translation length of its
generator.
If it is loxodromic let $\Gamma$ be its axis and let 
$h(\gamma)^*$ denote the quotient of $\Gamma$ in $N$, 
and if it is parabolic or has translation length
less than $\ep_0$, let $\til\T_h(\ep)$ (for $0<\ep\le\ep_0$) be the
invariant lift of its corresponding $\ep$-Margulis tube in $N$.

\begin{lemma}{Curve shortening}
Suppose $h:\gamma\to N$ is as above and satisfies the condition
$$
l(\gamma) \le \ell_N(h) + C.
$$
Let $\alpha\subset \til \gamma$ be an arc whose length is
$l(\alpha) = n l(\gamma) + r$, with $n\in\Z_{\ge 0}$ and
$r\in[0,l(\gamma))$.
Let $\ep\in(0,\ep_0]$ be given.
Then $h(\alpha)$ can be deformed, rel endpoints, to a path
$p_0 * \alpha' * p_1$ with the following properties:
\begin{enumerate}
\item $p_i$ are geodesics orthogonal to $\Gamma$ or $\til \T_h(\ep)$, and
$l(p_i) \le a_1$, for $i=0,1$.
\item $\alpha'$ lies on the geodesic $\Gamma$ if 
$\ell_N(h) \ge \ep$, and in $\til\T_h(\ep)$ otherwise.
\item $l_N(\alpha') = n l' + r'$, where 
  \begin{equation}
    \label{l' is max}
            l'=\max(\ell_N(h), \ep)
  \end{equation}
and
\begin{equation}
  \label{r' bound}
|r-r'| \le a_2.
\end{equation}
\end{enumerate}
The constants $a_i$ depend only on $\ep$ and $C$.
\end{lemma}

\begin{proof}
Fix a generator $\bar\gamma$ of $\pi_1(\gamma)$, and let $A=h_*(\bar\gamma)$.
Let $\ell = \ell_N(h)$ be the translation distance of $A$. 
Suppose first that $A$ is loxodromic.
If $x\in\Hyp^3$ with $d(x,\Gamma) = r$ and $s=d(x,Ax)$ then one can
obtain (see e.g. Buser \cite[2.3.1]{buser:surfaces}
\begin{equation}
  \label{s lower}
  \sinh\frac{s}{2} \ge \sinh \frac{\ell}{2}\cosh r
\end{equation}
Suppose that $s\le \ell+C$ and $\ell\ge \ep$. Then, since the function
$\sinh(t+C)/\sinh t$ is decreasing, we have
$\cosh r \le \sinh(\ep+C)/\sinh(\ep)$, yielding an upper bound on $r$.

Let $\xi\in\Gamma$ be the closest point to $x$ (so
$[x,\xi]\perp\Gamma$), and let $x_t$ be the point on $[x,\xi]$ with
$d(x_t,\xi)=t$. Suppose that $\ell \le \ep$ and hence $s\le \ep+C$.
We claim that
\begin{equation}\label{xt upper}
d(x_t,A(x_t)) \le a\frac{\cosh t}{\cosh r}
\end{equation}
where $a$ depends only on $\ep+C$.
To see this, note that the distance from $x$ to $Ax$ along the
$r$-equidistant surface to $\Gamma$ is at most a constant $c'$, since we can
project the geodesic $[x,Ax]$ to this surface with bounded distortion.
On the other hand the distance from $x_t$ to $Ax_t$ along the
$t$-equidistant surface is  $\sqrt{\ell^2\cosh^2 t +
  \theta^2\sinh^2 t}$, where $\ell+i\theta$ is the complex translation
length of $A$. A little algebra yields (\ref{xt upper}).

It follows that there exists $a_1$ depending only on $\ep$ such that
$x_t\in \til T(\ep)$ for all $t\in[0,r-a_1]$.

When $A$ is parabolic the same discussion works using coordinates
adapted to its parabolic fixed point, replacing equidistant surfaces
to $\Gamma$ with horospheres.

\medskip

Now for any arc $\alpha \subset \til\gamma$, with endpoints $x,y$,
let $p_0$ and $p_1$ be the orthogonal geodesics from  $\til
h(x)$ and $\til h(y)$, respectively, to $\Gamma$ when $\ell \ge \ep$,
and to $\til T(\ep)$ when $\ell < \ep$. The above discussion bounds
the lengths of $p_0$ and $p_1$ uniformly. Let $\xi,\eta$ be the
other endpoints of $p_0$ and $p_1$.

\medskip

Consider first the case that $\ell\ge \ep$.
Let $\alpha'$ be the geodesic segment $[\xi,\eta]$ 
on $\Gamma$. Write $l(\alpha) = nl(\gamma)  + r$ as in the statement
of the lemma. 
If $r=0$ then $\eta = A^{\pm n}\xi$, and it is immediate
that $l(\alpha') = n\ell$, which is what we wished to prove.

Suppose $n=0$, so that
$l(\alpha) = r < l(\gamma)$, and let us prove (\ref{r' bound}). 
Identify $\til\gamma$ with $\R$ so that the generator $\bar\gamma$ of
$\pi_1(\gamma)$ acts by $t\mapsto t+l(\gamma)$.
Let $x'= x+l(\gamma)$. We may assume $y\in[x,x']$. Let $\xi'\in
\Gamma$ be the closest point to $\til h(x')$.  Then
$\til h(x),\til h(x') $ and $\til h (y)$ are distance at most $a_1$
from $\xi,\xi'$ and $\eta$, respectively.

We have 
\begin{equation}\label{xi upper}
d(\xi,\eta) \le 2a_1 + |x-y|
\end{equation}
and
\begin{equation}\label{xi' upper}
d(\xi',\eta) \le 2a_1 + |x'-y|
\end{equation}
by the triangle inequality and the fact that $\til h$ is arclength-preserving.
We have by hypothesis
\begin{equation}\label{total lengths close}
|x-x'| - d(\xi,\xi') \le C.
\end{equation}
We have $|x'-y| = |x'-x|-|y-x|$ since
$y\in[x,x']$, and we have 
$d(\xi',\eta) \ge d(\xi,\xi') - d(\xi,\eta)$ by the triangle
inequality, so (\ref{xi' upper}) yields:
\begin{equation}
  \label{xi lower pre}
  d(\xi,\xi') - d(\xi,\eta) \le 2a_1 + |x'-x| - |y-x|.
\end{equation}
Now rearranging and applying (\ref{total lengths close}) we have
\begin{equation}
  \label{xi lower}
  d(\xi,\eta) \ge  |y-x| - 2a_1 - C.
\end{equation}
Now $d(\xi,\eta)=r'$ and $|y-x|=r$, so 
(\ref{xi upper}) and (\ref{xi lower}) together 
bound $|r-r'|$. 

Together with the upper bound (\ref{xi' upper}) on $d(\xi',\eta)$, we
find that $\eta$ is in fact 
constrained to within bounded distance of the point in
$[\xi,\xi']$ at distance $r$ from $\xi$ (as opposed to the point on
the other side of $\xi$ at the same distance). This together with the 
case of $n>0,r=0$ yields the general case, by concatentation.

\medskip
Now consider the case that $\ell<\ep$. 
Writing $l(\alpha) = nl(\gamma) + r$, let $\alpha'$ be the
concatenated chain of geodesic segments connecting the sequence of
$\xi=A^0\xi,\ldots,A^n\xi,\eta$. The length estimate on $\alpha'$ is immediate.

\end{proof}

\segment{Truncated curves}

\begin{lemma}{short truncation}
Let $\gamma$ be a simple closed geodesic 
or a simple properly embedded geodesic line
on a hyperbolic surface $S$.
Suppose that $\gamma$ avoids the $\ep_1$-thin part of
$S$, except possibly for its ends leaving the $\ep_1$-cusps.
Given
$\ep$, there exists a homotopically non-trivial, 
non-peripheral curve $\gamma_\ep$, homotopic to a simple curve,  with
the following properties:  
\begin{enumerate}
\item $\gamma_\ep$ has length at most $L_2$.
\item $\gamma_\ep$ is composed of at most 2 segments of $\gamma$ and
  at most 2 bridge arcs for $\gamma$ of length at most $\ep$.
\end{enumerate}
The constant $L_2$ depends only on $\ep,\ep_1$ and the topology of
$S$. 
\end{lemma}

\begin{pf} (Sketch)
We may assume $\ep<<\ep_1$. 
Let $L = 3\pi|\chi(S)|/\ep$, and let $L_2 = 2(L+\ep)$.

First, in case $\gamma$ is already a simple closed geodesic with
length less than $L_2$, take $\gamma_\ep = \gamma$.

If $\gamma$ is a properly embedded geodesic whose length outside the
$\ep$-cusps is at most $L$, $\gamma_\ep$ can be
constructed as a boundary component of a slight thickening of $\gamma$ union
the cusps (see Figure \ref{truncation cusp case}).

\realfig{truncation cusp case}{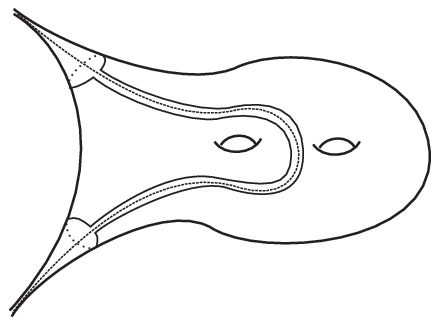}{Constructing $\gamma_\ep$
  when $\gamma$ is a bounded arc modulo cusps.}

Finally if $\gamma$ has length larger than $L$ outside the
$\ep$-cusps, let $\gamma_0$ denote any arc of $\gamma$
outside the
cusps of length $L$, and consider the immersion of
$\gamma_0\times[-\ep,\ep]$ taking  
each $\{p\}\times[-\ep,\ep]$ to the geodesic segment of
length $2\ep$ orthogonal to $\gamma$ at its midpoint $p$. The choice
of $L_2$ and the bound of $2\pi|\chi(S)|$  on the area of $S$ imply
that there must be a point $z$ covered with multiplicity 3 by this map. 
Thus $\gamma$ passes within $\ep$ of $z$ three times with (unoriented)
directions
varying by $O(\ep)$ in the projectivized tangent bundle. At least two
of those times the orientations 
match, and we may build $\gamma_\ep$ as shown in Figure \ref{two
  surgeries}.  The resulting curve is homotopically nontrivial and
non-peripheral because it can be smoothed out to have curvature near 0
and length bounded away from 0.
\end{pf}
\realfig{two surgeries}{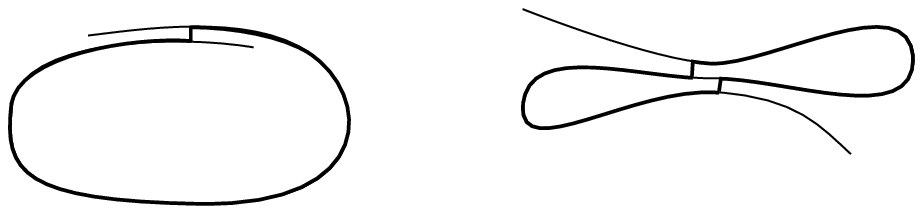}{Constructing $\gamma_\ep$ when
$\gamma$ is long}

\segment{Junctures of convex sets}

\begin{lemma}{convex juncture}
Fix $r_0>0$. For any $b>0$ there is a $d>0$ such that the following
holds: 
if $A$ and $B$ are convex sets in $\Hyp^3$, and suppose that
$\dist(A,B) \ge r_0$. Then
$$
\diam(\NN_b(A) \intersect \NN_b(B)) \le d.
$$
\end{lemma}
Here $\NN_b$ denotes $b$-neighborhood in $\Hyp^3$.

\begin{pf}
Convexity and the lower bound on $\dist(A,B)$ imply that $A$ and $B$
are contained in disjoint half-spaces $H_A$ and $A_B$ that are at
least $r_0$ apart.  The intersection $\NN_b(H_A)\intersect \NN_b(H_B)$
is compact so it has some diameter bound $d$. A sharp value is not
hard to compute but we will not need it.
\end{pf}

\providecommand{\bysame}{\leavevmode\hbox to3em{\hrulefill}\thinspace}

\end{document}